
%
%
\documentstyle{amsppt}
\TagsOnRight
%
\vsize=225truemm
\hsize=138truemm
\parindent=12pt 		  
\parskip=0pt
%
%
%
%
\loadbold
\loadeusm
\def\scr#1{{\fam\eusmfam\relax#1}}
\NoBlackBoxes
%
\def\lines#1:{\ifhmode\vskip#1\baselineskip\noindent\leftskip=0pt
              \else\ifvmode\vskip#1\baselineskip
                    \fi
              \fi}
\def\openC{\Bbb C}
\def\openG{\Bbb G}

\def\openP{\Bbb P}

\def\openZ{\Bbb Z}
\def\Hom{\operatorname{Hom}\,}

\def\coker{\operatorname{coker}\,}
\def\rk{\operatorname{rk}\,}
\def\Im{\operatorname{Im}\,}
\def\dim{\operatorname{dim}\,}
\def\ker{\operatorname{ker}\,}

\def\deg{\operatorname{deg}\,}
\def\Sing{\operatorname{Sing}\,}

\def\T{\operatorname{T}}
%
%
%
\topmatter
\title
On surfaces in ${\openP}^4$  and 3-folds in ${\openP}^5$
\endtitle
%
\author
Wolfram Decker, Sorin Popescu
\endauthor

\address
Wolfram Decker\hfil\break
Fachbereich Mathematik,
Universit\"at des Saarlandes,
D 66041 Saarbr\"ucken,
Germany
\endaddress
\email
decker\@math.uni-sb.de
\endemail
\address
{Sorin Popescu\hfil\break
Fachbereich Mathematik,
Universit\"at des Saarlandes,
D 66041 Saarbr\"ucken,
Germany}
\endaddress
\email
popescu\@math.uni-sb.de
\endemail
%

%
\toc
\widestnumber\head{2.}
\head 0. Introduction\endhead
\head 1. Constructions via syzygies\endhead
\head 2. Liaison\endhead
\head 3. Adjunction theory\endhead
\head 4. Surfaces in $\openP^4$\endhead
\head 5. 3-folds in $\openP^5$\endhead
\head 6. Examples: Two families of birational Calabi-Yau 3-folds
in $\openP^5$\endhead
\head 7. Overview\endhead
\endtoc
%
\endtopmatter
\document
\baselineskip=15pt 		
%
\head{0. Introduction}
\endhead
\noindent We report on some recent progress in the classification of
smooth projective varieties with small invariants. This progress is
mainly due to the finer study of the adjunction mapping by Reider,
Sommese and Van de Ven \cite{So1}, \cite {VdV}, \cite{Rei}, \cite{SV}.
Adjunction theory is a powerful tool for determining the type of a
given variety.
Classically, the adjunction process was introduced by Castelnuovo and
Enriques \cite{CE} to study curves on ruled surfaces. The italian geometers
around the turn of the century also started the classification of smooth
surfaces in $\openP^4$ of low degree. Further classification results are due
to Roth \cite{Ro1}, who uses the adjunction mapping to get surfaces with
smaller invariants already known to him (compare \cite{Ro2} for
adjunction theory on 3-folds). Nowadays, through the effort of
several mathematicians, a complete classification of smooth surfaces in
$\openP^4$ and smooth $3$-folds in $\openP^5$ has been worked out up
to degree $10$ and $11$ resp. Moreover, in the $3$-fold case the
classification is almost complete in degree $12$.
For references see section 7.\par
\noindent One motivation to study these varieties comes from
Hartshorne's conjecture \cite{Ha1}. In the case of codimension $2$
this suggests that already smooth $4$-folds in $\openP^6$ should be
complete intersections. Another motivation originates from two
mutually corresponding finiteness results. Ellingsrud and Peskine
\cite{EP} proved that there are only finitely many families of smooth
surfaces in $\openP^4$ which are not of general type. However, the
question of an exact degree bound $d_0$ is still open. By \cite{BF} $d_0
\leq 105$. Examples are known only up to degree $15$ and one
actually believes that $d_0 = 15$. The analogous finiteness result
holds for $3$-folds in $\openP^5$ \cite{BOSS1}. In this case one
expects a much higher degree bound. Nevertheless examples had been
known so far only up to degree $14$ \cite{Ch3}. In this note we
present, among other things, three new smooth $3$-folds in $\openP^5$
of degree $13, 17$ and $18$ resp.\par
\noindent How to construct examples ?\par
\noindent Let us recall that every smooth projective variety of
dimension $m$ can be embedded in $\openP^{2m+1}$. So e.g. in the
surface case we could try to work with general projections from
points in $\openP^5$. However Severi's theorem \cite{Se} tells us that
every non-degenerate smooth surface in $\openP^4$ except the Veronese
surface is linearly normal. Similarly by Zak's theorem \cite{Za} every
non-degenerate smooth $3$-fold in $\openP^5$ is linearly normal.\par
\noindent There are two other classical construction methods. One is
to study linear systems on abstract varieties. This works especially
well for rational, abelian and bielliptic surfaces. The other is
liaison \cite{PS} starting with a known local complete
intersection variety (presumably of lower degree). With
a few exceptions these methods failed to produce examples in higher
degree. In the case of liaison this is mainly due to the fact, that the
varieties to be constructed tend to be minimal in their even liaison class
(compare \cite{LR}). Whereas, if we consider e.g. linear systems of
curves on minimal surfaces, the base points have to be in a special
position. Such configurations are hard to find.\par
\noindent In this context a new construction method for surfaces
$X\subset \openP^4$ (more generally
$(n-2)$-folds $X\subset\openP^n$) was introduced in \cite{DES}
(compare also \cite {Po}). The basic
idea is an application of Beilinson's spectral sequence \cite{Bei}:
To construct the ideal sheaf $\scr J_X$ and thus $X$ one has to
construct the Hartshorne-Rao modules of $X$
first. Involving corresponding syzygy bundles as suggested by the
spectral sequence one finds vector bundles $\scr F$ and $\scr G$ on
$\openP^n$ with $\rk \scr G = \rk \scr F + 1$, and a morphism
$\varphi\in\Hom (\scr F, \scr G)$, whose minors define the desired $X$.
{}From the syzygies of the Hartshorne-Rao modules one can compute the
syzygies of $\scr J_X$ and so the explicit equations
for $X$. Typically, part of the geometry behind $X$ can already be seen
from the syzygies. The smoothness of $X$ can be checked via the implicit
function theorem, i.e., by a straightforward computation. Since these
computations are very extensive one has to rely on a computer and a computer
algebra system. Currently, Macaulay \cite{Mac} is the only system which is
powerful enough to handle the computations.\par
\noindent In some cases $X$ is not minimal in its even liaison class,
or a minimal element in the complementary even liaison class has low degree
and can be identified. In fact, by studying the equations we find examples
where $X$ can be reconstructed via liaison from a reducible
scheme $X'$ of lower degree.
It is hard to find such reducible schemes a priori.
\proclaim {Problem} Find a geometric construction for all
examples constructed via syzygies.\qquad\qed\endproclaim
\noindent
{\bf Notation.}\quad $R =\openC [x_0,\dots , x_n] = \underset
{q\in\openZ}\to{\textstyle\bigoplus} S^qV^*$ will be
the homogeneous coordinate ring of $\openP^n$, so
$H^0(\scr O_{\openP^n}(1))=V^*$.
If $X\subset \openP^n$ is a fixed smooth subvariety, then $d$ will denote its
degree, $\pi$ its sectional genus, $H$ the hyperplane class and
$K$ the canonical class.\qquad\qed
\smallskip\noindent
{\bf Acknowledgements.}\quad Both authors are grateful to Frank-Olaf Schreyer
for many helpful conversations. We also thank Mark Gross for interesting
discussions. The first author, who lectured in Tokyo and Durham on the topic of
this note, would like to thank the Japan Society for the
Promotion of Science and the Deutscher Akademischer Austauschdienst
for their support and Masaki Maruyama and the University of Kyoto for
their hospitality. Last but not least we would like to thank the
organizers of the Durham symposium on Vector Bundles in Algebraic Geometry
for creating a stimulating atmosphere during this meeting.\qquad\qed
\head{1. Constructions via syzygies}
\endhead
\noindent Following \cite{DES} we want to construct a codimension 2
subvariety $X\subset \openP^n$ as the determinantal locus of a map
between vector bundles. So we are looking for vector bundles $\scr F$
and $\scr G$ on $\openP^n$ with $\rk \scr F = f$ and $\rk \scr G =
f+1$, and a morphism $\varphi \in\Hom (\scr F, \scr G)$ whose minors
vanish in the expected codimension 2. In this case $X=V(\varphi)$ is
a locally Cohen-Macaulay subscheme and the Eagon-Northcott complex \cite{BE}
$$
0\leftarrow \scr O_X (m)\leftarrow \scr O (m)\cong \bigwedge^f
\scr F^*\otimes \bigwedge^{f+1} \scr G\leftarrow \scr G
\overset\varphi\to\leftarrow \scr F\leftarrow 0
$$
is exact and identifies $\coker \varphi$ with the twisted ideal sheaf
$$
\coker \varphi\cong \scr J_X (m),\qquad m = c_1\scr G - c_1\scr F.
$$
Furthermore, a mapping cone between the minimal free
resolutions of $\scr F$ and $\scr G$ is a (not necessarily minimal)
free resolution of $\scr J_X(m)$. So for a given $\varphi$ we can
derive an explicit system of homogeneous equations for its dependency
locus $X$.
\medskip\noindent
{\bf Remark 1.1.}\quad Let $\varphi_1,\varphi_2\in\Hom (\scr F, \scr
G)$ be morphisms whose minors vanish in codimension 2. Then
$V(\varphi_1)$ and $V(\varphi_2)$ lie in the same irreducible
component of the Hilbert scheme (compare e.g.
\cite{BB}, \cite{MDP}).\qquad\qed\par\medskip
\noindent To construct a variety with the desired numerical
invariants one has to find appropriate $\scr F$ and $\scr G$. Clearly
$\scr F$ and $\scr G$ reflect the structures
of the graded finite length $R$-modules
$$
H_*^i \scr J_X = \underset {q\in\openZ}\to{\textstyle\bigoplus} H^i
(\openP^n, \scr J_X(q)),\qquad i=1,\dots,\dim X,
$$
called the Hartshorne-Rao modules of $X$.
E.g., $X$ is  projectively Cohen-Macaulay, i.e., its Hartshorne-Rao
modules are trivial, iff $\scr F$ and $\scr G$ can be chosen to be
direct sums of line bundles. Or compare \cite{Ch2} for the
$\Omega$-resolution of a projectively Buchsbaum variety. In this case, in
particular,
the multiplication maps of the Hartshorne-Rao modules are trivial.
\medskip\noindent
{\bf Remark 1.2.}\quad  Smooth projectively Cohen-Macaulay and smooth
projectively Buchsbaum varieties of codimension 2, which are not of general
type, are completely classified (compare \cite{Ch3}).\qquad\qed\medskip
\noindent In any case it is a natural idea to construct the Hartshorne-Rao
modules first. Then one may involve corresponding syzygy bundles as direct
summands in order to find $\scr F$ and $\scr G$. Recall:\par
\proclaim {Proposition 1.3} Let $M=\underset{q\in \openZ}\to\bigoplus M_q$ be a
graded $R$-module of finite length and let
$$
0\leftarrow M\leftarrow L_0\overset{\alpha_1}\to\leftarrow
L_1\leftarrow\dots \overset{\alpha_{n+1}}\to\leftarrow L_{n+1}\leftarrow 0
$$
be its minimal free resolution. Then, for $1\leq i\leq n-1$, the
sheafified syzygy module
$$
\scr F_i = \scr Syz_i (M) = (\ker \alpha_i)^\sim = (\Im \alpha_{i+1})^\sim
$$
is a vector bundle on $\openP^n$ with the intermediate cohomology
$$
\underset {q\in\openZ}\to{\textstyle\bigoplus} H^j \left(\openP^n, \scr F_i
(q)\right) = \left\{
\aligned
M&\quad j=i\\
0&\quad j\neq i,\quad 1\leq j\leq n-1\quad .
\endaligned
\right .
$$
Conversely, any vector bundle $\scr F$ on $\openP^n$ with this
intermediate cohomology is stably\par
\noindent equivalent with $\scr F_i$, i.e.,
$$
\scr F \cong \scr F_i\oplus \scr L\ ,\qquad \scr L \text{\quad a
direct sum of line bundles.}\qquad \qed
$$
\endproclaim
\noindent
{\bf Example 1.4.}\quad Consider $\openC$ as a graded $R$-module
sitting in degree 0. The minimal free resolution of $\openC\,(i)$ is
the Koszul complex
$$
0\leftarrow \openC\,(i)\leftarrow \bigwedge^0 V^*{\textstyle\otimes}
R(i) \leftarrow \dots \leftarrow\bigwedge^{n+1}
V^*{\textstyle\otimes} R(i-n-1)\leftarrow 0\ .
$$
The corresponding syzygy bundles are the twisted bundles of
differentials,
$
\scr Syz_i (\openC\,(i)) \cong \Omega^i(i)\ .
$
It follows from the sheafified Koszul complex, that
$
\Hom (\Omega^i(i), \Omega^j(j)) \cong \bigwedge^{i-j} V\ ,
$
the isomorphisms being given by contraction (cf. \cite{Bei}).\qquad\qed\par
\medskip\noindent Which syzygy bundles should be involved in the construction
of $\scr F$ and $\scr G$ ? This can be found out by analyzing
Beilinson's spectral sequence for $\scr J_X (m)$. Recall:
\medskip
\proclaim {Theorem 1.5} \cite{Bei}\quad
For any coherent sheaf $\scr S$ on $\openP^n$ there is a spectral
sequence with $E_1$-terms
$$
E_1^{pq} = H^q\left(\openP^n, \scr S(p)\right)\otimes \Omega^{-p}(-p)
$$
converging to $\scr S$, i.e., $E_\infty^{pq} = 0$ for $p+q\neq 0$ and
$\bigoplus E_\infty^{-p,p}$ is the associated graded sheaf of a suitable
filtration of $\scr S$.\qquad\qed
\endproclaim
\noindent This theorem is often used to construct $\scr S$ by determining the
differentials of the spectral sequence first. A crucial point is
that the $d_1$-differentials
$$
\aligned
d_1^{pq} &\in \Hom (H^q(\openP^n, \scr S(p))\otimes \Omega^{-p}(-p),\
H^q(\openP^n, \scr S(p+1))\otimes \Omega^{-p-1}(-p-1))\\
&\cong \Hom (V^* \otimes
H^q(\openP^n, \scr S(p)),\  H^q(\openP^n, \scr S(p+1))
\endaligned
$$
coincide with the natural multiplication maps.
In our case $\scr S = \scr J_X(m)$, and we will interpret
one part of Beilinson's spectral sequence
as the spectral sequence of a vector bundle $\scr F$, the other part
as that of a vector bundle $\scr G$. The differential between the two
parts will give the morphism $\varphi : \scr F\rightarrow \scr G$
whose cokernel is the desired $\scr J_X(m)$. The twist m will be mainly
$n$ or $n-1$ (compare \cite{DES, 1.7} for the corresponding Beilinson
cohomology tables in the surface case).
\medskip
\noindent How to check the smoothness of $X$ ? If the bundle $\Hom
(\scr F, \scr G)$ is globally generated and $n\le 5$, then we know from
\cite{Klm}, that the generic
$\varphi\in\Hom (\scr F, \scr G)$ gives rise to a smooth $X$.
This works well, if $X$ is projectively Cohen-Macaulay.
Similarly, if $X$ is projectively Buchsbaum, we may apply \cite{Ch1}.
In the general case however, we mostly have to rely on a computer as
explained in the introduction.
\medskip\noindent
{\bf Example 1.6.}\quad We will construct a family of smooth 3-folds $X\subset
\openP^5$ with the numerical invariants $d=18$, $\pi = 35$,
$\chi(\scr O_X)=2$ and $\chi(\scr O_S)=26$, where $S$ is a general
hyperplane section of $X$. Let us analyze Beilinson's spectral sequence for
$\scr J_X (4)$.
We first need information on the dimensions $h^i \scr J_X (m)$, $m =
-1,\dots,4$.
In view of Riemann-Roch a plausible Beilinson cohomology table is
%
$$
\vbox{\offinterlineskip
\halign{
\vrule height10.5pt depth 5.5pt#
&\hbox to 25pt{\hfil$#$\hfil}
&\vrule#
&\hbox to 25pt{\hfil$#$\hfil}
&\vrule#
&\hbox to 25pt{\hfil$#$\hfil}
&\vrule#
&\hbox to 25pt{\hfil$#$\hfil}
&\vrule#
&\hbox to 25pt{\hfil$#$\hfil}
&\vrule#
&\hbox to 25pt{\hfil$#$\hfil}
&\vrule#
&\hbox to 30pt{\hfil}#
&#
\cr
\omit&&\omit\hbox to
0pt{\hskip-3pt\hbox{$\bigg\uparrow$}\raise5pt\hbox{$i$}\hss}\cr
\multispan{13}\hrulefill\cr
& && && && && && && \cr
\multispan{13}\hrulefill\cr
&24&& && && && && &\cr
\multispan{13}\hrulefill\cr
& &&1&&6&&3&& && &\cr
\multispan{13}\hrulefill\cr
& && && && && && &\cr
\multispan{13}\hrulefill\cr
& && && && && && &\cr
\multispan{13}\hrulefill\cr
& && && && && && &\cr
\multispan{14}\hrulefill
&\vbox to0pt{\vss\vskip5.5pt\hbox{$\!-\negthickspace\negmedspace@>> m
>$}\vss}\cr
}}
$$
%
\noindent Suppose that a smooth 3-fold $X$ with these data exists. Then
Beilinson's theorem yields an exact sequence
$$
0\rightarrow \scr F = 24 \scr O (-1)\rightarrow \scr G \rightarrow
\scr J_X (4) \rightarrow 0\ ,
$$
where $\scr G$ is the cohomology of a monad
$$
0\rightarrow \Omega^4(4) \overset {d_1^{-4,3}}\to\rightarrow
6\Omega^3(3) \overset{d_1^{-3,3}}\to\rightarrow 3
\Omega^2(2) \rightarrow 0\ .
$$
On the other hand, the generic module with Hilbert function $(1,6,3)$
has syzygies of type\par
%
%
$$
\vbox{%
\halign{&\hfil\,$#$\,\hfil\cr
0\leftarrow M\leftarrow R(4)\ \cr
&\vbox to 10pt{\vskip-4pt\hbox{$\nwarrow$}\vss}\ 18 R(2)\leftarrow 52
R(1)\leftarrow 60 R&&24 R(-1)\cr
&&\vbox to 10pt{\vskip-4pt\hbox{$\nwarrow$}\vss}&\oplus\cr
&&&10 R(-2)&&\vbox to 10pt{\vskip-4pt\hbox{$\nwarrow$}\vss}&
12 R(-3)\leftarrow 3 R(-4)\leftarrow 0\ .\cr
}}
$$
%
A check on the ranks and the intermediate cohomology of
$\scr G$ and $\scr Syz_3 (M)$ suggests that conversely
it is promising to start with $\scr F = 24 \scr O (-1)$ and
$\scr G = \scr Syz_3 (M)$.
Indeed, for the map $\varphi \in \Hom (24 \scr O (-1), \scr
Syz_3 (M))$ given by the syzygies,  one can check that the
minors of $\varphi$ vanish along a smooth 3-fold $X$.
By construction $\scr J_X$ has syzygies of type
$$
0\leftarrow \scr J_X \leftarrow 10 \scr O(-6)\leftarrow 12 \scr O
(-7) \leftarrow 3 \scr O (-8) \leftarrow 0\ .
$$
In particular, $X$ is cut out by 10 sextics. From the syzygies it follows
that $\omega_{X}(1) = \scr Ext^2 (\scr O_{X}, \scr O(-6))(1)$ is a quotient
of ${24\scr O}$, and since $(K+H)^2\cdot K=-4$ (compare section 5) we deduce
that the Kodaira dimension $\kappa(X)=-\infty$.
\qquad\qed
\head{2. Liaison}
\endhead
\noindent We recall the definition and some basic results \cite{PS}.
Let $X,X'\subset\openP^n$ be two locally Cohen-Macaulay subschemes of
pure codimension $2$ with no irreducible components in common. $X$ and
$X'$ are said to be (geometrically) {\it linked} $(r,s)$, if there exist
hypersurfaces $V_1$ and $V_2$ of degrees $r$ and $s$ resp. such that
$
X\cup X' = V_1\cap V_2\ .
$
Then there are the standard exact sequences
$$
0 \rightarrow \omega_X\rightarrow \scr O_{V_1\cap V_2}(r+s-n-1)
\rightarrow \scr O_{X'}(r+s-n-1)\rightarrow 0\ ,
$$
$$
0 \rightarrow \omega_X\rightarrow \scr O_X(r+s-n-1)
\rightarrow \scr O_{X\cap X'}(r+s-n-1)\rightarrow 0\ .
$$
The degrees and sectional genera of $X$ and $X'$ are related by
$$
d+d' = r\cdot s \text{\qquad and\qquad} \pi -\pi' = {1\over 2} (r+s-4)(d-d')\ ,
$$
and
$$
\chi(\scr O_{X'}) = \chi (\scr O_{V_1\cap V_2}) - \chi (\scr
O_X(r+s-n-1))\ .
$$
Under suitable assumptions (e.g., if $H^1(\scr F(r))=H^1(\scr F(s))=0$)
we may deduce from a given locally free resolution
$
0 \rightarrow \scr F\rightarrow \scr G\rightarrow \scr J_X\rightarrow
0
$
of $\scr J_X$ a resolution
$$
0\rightarrow\scr G\,\check{}\,(-r-s)\rightarrow \scr
F\,\check{}\,(-r-s)\oplus \scr O(-r)\oplus\scr O(-s)\rightarrow\scr
J_{X'} \rightarrow 0
$$
of $\scr J_{X'}$ by taking a mapping cone as in \cite{PS, Prop. 2.5}.
Moreover, the Hartshorne-Rao modules of $\scr J_{X'}$ are
$\openC$-dual to those of $X$:
$$
H_*^{n-1-i}\scr J_X \cong (H_*^i\scr J_{X'})^* (n+1-r-s),\qquad
i=1,\dots, n-2\ .
$$
Liaison can be used to construct new subvarieties starting from given
ones. Hence it is useful to know, under which conditions a residual
intersection will be smooth. One result in this direction is a
special case of \cite{PS, Prop.4.1}:
\proclaim {Theorem 2.1} (Peskine-Szpiro). Let $X\subset \openP^n$,
$n\leq 5$, be a local complete intersection of codimension $2$. Let
$m$ be a twist such that $\scr J_X (m)$ is globally generated. Then
for every pair $d_1,d_2\geq m$ there exist forms $f_i\in H^0(\scr
J_X(d_i))$, $i=1,2$, such that the corresponding hypersurfaces $V_1$
and $V_2$ intersect properly, $V_1\cap V_2 = X\cup X'$, where
\roster
\item"{(i)\ }" $X'$ is a local complete intersection,
\item"{(ii)\ }" $X$ and $X'$ have no common component,
\item"{(iii)\ }" $X'$ is nonsingular outside a set of positive
codimension in $\Sing X$.\qquad\qed
\endroster
\endproclaim
\head{3. Adjunction theory}
\endhead
\noindent In this section $(X,H)$ will denote a polarized pair, where
$X$ is a smooth, connected, projective variety of dimension $m\ge 2$ and
$H$ is a very ample divisor on $X$.\ \ $K=K_X$ will be a canonical
divisor on $X$. Before reviewing the general theory behind the
{\it adjunction map} $\Phi = \Phi_{\mid K+(m-1)H\mid}$, we will
give an example.
\medskip\noindent
{\bf Example 3.1.} \cite {Roo}\quad
Let
$
\varphi = (\varphi_{ij})_{\vbox{
\baselineskip 5.5pt
\vskip 4pt
\hbox{$\scriptstyle 0\leq i\leq 2$}
\hbox{$\scriptstyle 0\leq j\leq 3$}
}}
$
be a general $3\times 4$-matrix with linear entries in $\openC
[x_0,\dots, x_4]$. Then $X = V(\varphi)$ is a smooth surface
$X\subset\openP^4$ with $d=6$ and $\pi = 3$. Let $H$ be the
hyperplane class of $X$. By dualizing $\varphi$, we obtain the resolution
$$
0\leftarrow \omega_X(1)\leftarrow 3\scr O
\overset{{}^t\varphi}\to\leftarrow 4\scr O(-1)\leftarrow
\scr O(-4)\leftarrow 0\ .
$$
\noindent
So $|K+H|$ is base point free, $N = \dim |K+H| = 2$ and we
have a well-defined adjunction map
$\Phi = \Phi_{\mid K+H\mid}: X\rightarrow \openP^2\ .$
Let $y_0,y_1,y_2$ be coordinates on $\openP^2$. Then  graph
$(\Phi)\subset \openP^4\times \openP^2$ is given by the equations
$$
y_0\varphi_{0j}(x) + y_1\varphi_{1j}(x) + y_2\varphi_{2j}(x)=0
\,,\quad j = 0,\dots,3\ .
$$
We may rewrite these equations as
$$
x_0\psi_{j0}(y)+\dots+ x_4\psi_{j4}(y) = 0\,,\quad j = 0,\dots,3\ ,
$$
where
$
\psi = (\psi_{jk})_
{\vbox{
\baselineskip 5.5pt
\vskip 4pt
\hbox{$\scriptstyle 0\leq j\leq 3$}
\hbox{$\scriptstyle 0\leq k\leq 4$}
}}
$
has linear
entries in $\openC [y_0,y_1,y_2]$. The general fibre of $\Phi$ is
defined by four independent linear forms in $\openC [x_0,\dots,x_4]$.
Hence $\Phi$ is birational with positive dimensional fibres precisely
in the points where $\psi$ drops rank:
$$
0\rightarrow 4\scr O_{\openP^2}(-5) \overset\psi\to\rightarrow
5\scr O_{\openP^2}(-4) \rightarrow \scr J_Z \rightarrow 0\ .
$$
So $\Phi : X\rightarrow \openP^2$ expresses $X$ as the blowing up of
$10$ points in $\openP^2$ and $X$ is embedded by the quartics through
these points, i.e., by the $4\times 4$-minors of $\psi$. In other words
$$
H \equiv 4L - \sum\limits_{i=1}^{10} E_i
$$
(with obvious notations) and $X$ is a Bordiga surface.\qquad\qed
\medskip\noindent
The first general result deals with the existence of the adjunction
map. It is a consequence of \cite{So1}, \cite{VdV}.
\proclaim {Theorem 3.2} Let $(X,H)$ and $K$ be as above. Then
$|K+(m-1)H|$ is base point free unless\par
\roster
\item"{(i)\ }" $(X,\scr O_X(H)) \cong (\openP^m, \scr O_{\openP^m} (1))$ or
$(\openP^2, \scr O_{\openP^2} (2))$,
\item"{(ii)\ }" $(X,\scr O_X(H)) \cong (Q, \scr O_Q (1))$, where $Q\subset
\openP^{m+1}$ is a smooth hyperquadric,
\item"{(iii)\ }" $(X,\scr O_X(H))$ is a scroll over a smooth
curve.\qquad\qed
\endroster
\endproclaim
\noindent
If $|K+(m-1)H|$ is base point free, then we denote by
$$
\vbox{
\halign{&$#$\hfil\cr
X&\overset\Phi\to\longrightarrow&\quad\openP^N\cr
r\searrow&&\nearrow s\cr
&X'\cr
}}$$
the Stein factorization of the adjunction map $\Phi$.
$X'$ is normal, $r$ is connected and $s$ is finite.

\proclaim{Theorem 3.3} \cite{So2}
Let $(X,H)$ and $K$ be as above and suppose that $|K+(m-1)H|$ is
base point free. Then there are the following possibilities:\par
\roster
\item"{(i)\ }" $\dim \Phi (X) =0$, and $K \equiv -(m-1)H$, i.e., $X$ is Fano of
index $(m-1)$,
\item"{(ii)\ }" $\dim \Phi (X)=1$, and the general fibre of $r$ is a
smooth quadric $Q$ such that $H$ induces $\scr O_Q (1)$,
\item"{(iii)\ }" $\dim \Phi (X) = 2<m$ and $r$ exhibits $X$ as a scroll
over a smooth surface,
\item"{(iv)\ }" $\dim \Phi (X) = m$\ .\qquad\qed
\endroster
\endproclaim
\noindent If $\dim \Phi (X) = m$ we write $L' = r_*(H)$ , $K' = K_{X'}$ and $H'
= K' + (m-1)L'$. The next result tells us, that in this case $r$
contracts precisely the linear $\openP^{m-1}\subset X$ with normal
bundle $\scr O_{\openP^{m-1}}(-1)$ (necessarily disjoint).
\proclaim
{Theorem 3.4} \cite{So1},\cite{So2}
Suppose that $\dim \Phi (X) = m$. Then $r : X\rightarrow X'$ is the
blowing up of a finite number of points on the smooth projective
variety $X'$. $L'$ and $H'$ are ample and
$$
r^* (H') \equiv K + (m-1)H\ .\qquad\qed
$$
\endproclaim
\noindent
In the above situation $(X', L')$ is called the {\it {first reduction}}
of $(X, H)$ \cite{So5}.
\smallskip
\noindent When is $s$ an embedding ? The answer is given by
\proclaim
{Theorem 3.5} \cite{SV}
Suppose that $\dim \Phi (X) = m$. Then $H'$ is very ample, unless $X$
is a surface and
\roster
\item"{(i)\ }" $X = \openP^2 (p_1,\dots, p_7)$ and $H \equiv 6L -
\sum\limits_{i=1}^7 2 E_i$\quad (the Geiser involution),
\item"{(ii)\ }" $X = \openP^2 (p_1,\dots, p_8)$ and $H \equiv 6L -
\sum\limits_{i=1}^7 2 E_i - E_8$\,,
\item"{(iii)\ }" $X = \openP^2 (p_1,\dots, p_8)$ and $H \equiv 9L -
\sum\limits_{i=1}^8 3 E_i$\quad (the Bertini involution),
\item"{(iv)\ }" $X = \openP (\scr E)$,\quad where $\scr E$ is an indecomposable
rank
2 bundle over an elliptic curve, and $H \equiv 3B$, where $B$ is a section
with $B^2 = 1$ on $X$\ .\qquad\qed
\endroster\endproclaim
\noindent For surfaces the {\it adjunction process}, i.e., the study of
$|K+H|$,\ $|K'+ H'|$ etc., will finally lead to a minimal model. For
3-folds $X\subset\openP^5$ the situation is quite different. In this
case it is often successful to study $|K+H|$ instead of $|K+2H|$. Compare
section 5 for details and applications of further general results of
adjunction theory.
\head{4. Surfaces in $\openP^4$}
\endhead
\noindent In this section $X$ will denote a smooth non-degenerate surface in
$\openP^4$ and
$d=H^2$ its degree, $\pi = {1\over 2} H \cdot (K+H)+1$ its sectional
genus and $\chi =\chi(\scr O_X)=1-q+p_g$ its Euler characteristic.\par\noindent
$K^2$ may be computed from the double point formula (cf. \cite{Ha2, Appendix A,
4.1.3.})
$$
d^2 - 10d - 5 H\cdot K -2K^2 +12 \chi = 0\ .
$$
In order to classify surfaces of a given degree, one first has
to work out a finite list of admissible numerical invariants. One may
apply Halphen's upper bound for $\pi$ \cite{GP} in connection with the lifting
theorem
of Roth \cite{Ro1, p.152} and the following classification results:
\proclaim
{Theorem 4.1} \cite{Ro1}, \cite{Au}.
Let $X$ be contained in a hyperquadric $V^2\subset \openP^4$. Then
$\pi = 1 + \left[ d(d-4)/4\right ]$ and $X$ is either the complete intersection
of
$V^2$ with another hypersurface, or $X$ is linked to a plane in the
complete intersection of $V^2$ with another hypersurface.\qquad\qed
\endproclaim
\proclaim
{Theorem 4.2} \cite{Ko}, \cite{Au}.
Let $X$ be contained in an irreducible cubic hypersurface $V^3\subset
\openP^4$. Then either $X$ is projectively Cohen-Macaulay and linked
on $V^3$ to an irreducible scheme of degree $\leq 3$, or $X$ is
linked on $V^3$ to a Veronese surface, or to a quintic elliptic
scroll.\qquad\qed
\endproclaim
\proclaim
{Corollary 4.3} If $X$ is contained in a cubic hypersurface and if
$d\geq 9$, then $X$ is of general type.\qquad\qed
\endproclaim
\noindent To derive a lower bound for $\pi$ and bounds for $\chi$ we may use
Severi's Theorem \cite{Se} together with Riemann-Roch, the Hodge
index theorem, the Enriques-Kodaira classification and adjunction
theory. In the context of section 3 we note:
\proclaim
{Theorem 4.4} \cite{Au}, \cite{La}.
If $X$ is a scroll, then $X$ is a rational cubic or an elliptic
quintic scroll.\qquad\qed
\endproclaim
\proclaim
{Theorem 4.5} \cite{BR}, \cite{ES}.
If $X$ is a conic bundle, then $X$ is a Del Pezzo surface of degree 4,
or a Castelnuovo surface.\qquad\qed
\endproclaim
\proclaim
{Corollary 4.6}
If $d\geq 6$, then the adjunction map $\Phi$ is defined and $(K+H)^2
> 0$, i.e., $\dim \Phi (X) = 2$.\qquad\qed
\endproclaim
\noindent
Once the numerical invariants are fixed, we use the information on the
dimensions
$h^i\scr J_X (m)$ provided by Riemann-Roch and \cite{DES, 1.7}.
In some cases more information on the dimensions and the structures of the
Hartshorne-Rao
modules may be obtained by studying the relations between the multisecants to
$X$,
the plane curves on $X$ and the syzygies of $\scr J_X$ (compare \cite{PR}).
This information is
helpful for construction and classification purposes.
In other cases one has to go through the adjunction process  to
analyze, how a given $X$ fits into the Enriques-Kodaira classification. In any
case it is
crucial to know the number of exceptional lines on $X$. Le Barz' 6-secant
formula \cite{LB}
tells us, that the number of 6-secants to $X$ (if finite) plus the number of
exceptional
lines equals a polynomial expression $N_6 = N_6 (d,\pi,\chi)$\ \ (if $X$ does
not contain a
line with self-intersection $\geq 0$). This fits well with the ideas of section
1.
Once having constructed $X$ explicitly , we can compute the 6-secants easily.
For examples we refer to \cite{DES},\cite{Po}.\par\medskip\noindent
With the following example we would like to demonstrate, that the construction
via
syzygies is not always as straightforward as in Example 1.6.
\medskip\noindent {\bf Example 4.9.} \cite{Po}\quad
Let us construct a family of smooth surfaces $X\subset \openP^4$ with
$d=11$, $\pi =11$ and $\chi = 3$. In view of \cite{DES, 1.7} a
plausible Beilinson cohomology table for $\scr J_X (4)$ is
%
%
$$
\vbox{\offinterlineskip
\halign{
\vrule height10.5pt depth 5.5pt#
&\hbox to 25pt{\hfil$#$\hfil}
&\vrule#
&\hbox to 25pt{\hfil$#$\hfil}
&\vrule#
&\hbox to 25pt{\hfil$#$\hfil}
&\vrule#
&\hbox to 25pt{\hfil$#$\hfil}
&\vrule#
&\hbox to 25pt{\hfil$#$\hfil}
&\vrule#
&\hbox to 30pt{\hfil}#
&#
\cr
\omit\hbox to 0pt{\hskip-3pt\hbox{$\bigg\uparrow$}\raise5pt\hbox{$i$}\hss}\cr
\multispan{11}\hrulefill\cr
& && && && && && \cr
\multispan{11}\hrulefill\cr
&2&& && && && & \cr
\multispan{11}\hrulefill\cr
& &&1&& && && & \cr
\multispan{11}\hrulefill\cr
& && &&1&&4&&3& \cr
\multispan{11}\hrulefill\cr
& && && && && & \cr
\multispan{12}\hrulefill
\vbox to0pt{\vss\vskip5.5pt\hbox{$\!-\negthickspace\negmedspace@>> m
>$}\vss}\cr
}}
$$
%
\medskip\noindent
Suppose that a smooth surface $X$ with these data exists. Then
Beilinson's theorem yields a resolution of type
$$
0\rightarrow \scr F = 2 \scr O (-1) \oplus \Omega^3 (3) \overset \varphi
\to \rightarrow \scr G\rightarrow \scr J_X (4)\rightarrow 0\ ,
$$
where $\scr G$ is the cohomology of a monad
$$
0\rightarrow \Omega^2 (2) \overset {d_1^{-2,1}}\to \rightarrow 4
\Omega^1 (1) \overset {d_1^{-1,1}}\to \rightarrow 3 \scr O \rightarrow
0\ .
$$
Arguing as in example 1.6, we conversely choose
$\scr G = \scr Syz_1 (M)$, where $M$ is a module with Hilbert
function $(1,4,3)$ and a minimal free presentation of type
$$
0 \leftarrow M \leftarrow S(2) \overset {(\alpha, \beta
)}\to\leftarrow S(1) \oplus 7 S\ .
$$
So $M$ is the tensor product of the Koszul complex given by the
linear form $\alpha$ and the module $M'$ presented by $\beta$. We may
assume that $\alpha = x_4$ and that $M'$ is a module over $R' =
\openC [x_0,\dots , x_3]$.\ $M'$ has the same Hilbert function as
$M$, namely $(1,4,3)$. The general such $M'$ has syzygies of type
%
%
$$
\vbox{%
\halign{&\hfil\,$#$\,\hfil\cr
0\leftarrow M'\leftarrow R'(2)\ \cr
&\vbox to 10pt{\vskip-4pt\hbox{$\nwarrow$ \raise 5pt\hbox to
0pt{\hss$\beta$}}\vss}\cr
&&7 R'&&8 R'(-1)&&aR'(-2)\cr
&&&\vbox to 10pt{\vskip 0pt\hbox{$\nwarrow$}\vss}&\oplus
&\longleftarrow&\oplus\cr
&&&&(3+a)R'(-2)&&8 R'(-3)&&\vbox to 10pt{\vskip-6pt\hbox{$\nwarrow$
\raise 5pt \hbox to 3pt{\hss$^t\gamma$}}\vss}&
3 R'(-4)\leftarrow 0\ ,\cr
}}
$$
%
with $a=0$. It is easy to see, that in this case no morphism $\varphi\in\Hom
(\scr F,\scr G)$ is injective. The trick for the construction of $X$
is to choose $\beta$ special in order to obtain some extra syzygies
and thus a larger space $\Hom (\scr F, \scr G)$. We will construct a
module $M'$ with the above type of syzygies and $a=1$.
Equivalently, we will construct the $\openC$-dual module ${M'}^*$ by
defining its presentation matrix $\gamma = (\gamma_1, \gamma_2)$. Choose
four general lines $L_1,\dots ,L_4$ in the hyperplane $V(x_4)$, denote by
$\epsilon$ the presentation matrix in the direct sum of the
four Koszul complexes built on these lines and let
$\delta$ be a general $3\times 4$-matrix with entries in $\openC$.
Then $\epsilon$ and thus also $\gamma_1 = \delta\epsilon$ has four linear
1-syzygies.
Let $\gamma_2$ be given by 3 general quadrics. Then $\gamma$ presents an
artinian module as desired. With these choices the generic
$\varphi \in\Hom (\scr F, \scr G)$ yields a smooth surface $X$ cut out
by 8 quintics and 4 sextics. In general it is a plausible guess and in many
cases true that the number of 6-secants to a surface in $\openP^4$ is precisely
the number of sextic generators of its homogeneous ideal. Indeed, in our case
it is easy to see that $L_1,\dots ,L_4$ are precisely the 6-secants
to $X$ \cite{Po, Proposition 3.32}. Le Barz' 6-secant formula gives
$N_6(11, 11, 3) = 5$. Hence there is one exceptional line on $X$. One
can show that there are no other exceptional curves
\cite{Po, Proposition 3.31}.
Since $K^2 = 1$ by the double point formula $X$ is  of general type.\qquad\qed
\medskip\noindent
In some cases it is quite subtle to construct artinian modules with the desired
graded
Betti numbers. From this point of view the most difficult surfaces are the
abelian and
bielliptic surfaces known so far \cite{ADHPR2}. These are also the surfaces
with the
most beautiful geometry behind (compare \cite{ADHPR1} and \cite{Hu2}).
The link between the geometry and the
syzygies is provided by the distribution of the $2$- and $3$-torsion points on
the
Heisenberg invariant elliptic normal curves in $\openP^4$. In turn, these
curves are
related to the Horrocks-Mumford bundle. Our knowledge on this bundle has
influenced
the construction of further families of surfaces
(compare \cite{ADHPR1, Thm. 32},
\cite{DES, 2.5}, \cite{Po, 4.1 and 7.4}).
\head{5. 3-folds in $\openP^5$}
\endhead
\noindent In this section $X$ will denote a smooth, non-degenerate 3-fold in
$\openP^5$, $S$ a general hyperplane section, $d = H^3$ its degree and
$\pi = {1\over 2} H^2\cdot (K+2H)+1$ its sectional genus.\newline
We have two double point formulae, one for $X$,
$$
K^3 = - 5d^2 + d (2\pi + 25) + 24 (\pi -1) - 36 \chi (\scr O_S) - 24
\chi (\scr O_X),
$$
and one for $S$, which may be rewritten as
$$
H\cdot K^2 = {1\over 2} d(d+1) - 9(\pi - 1) + 6 {\chi} (\scr O_S)
$$
(compare e.g. \cite{Ok2}). So the basic invariants of $X$ are
$d,\pi, \chi (\scr O_X)$ and $\chi (\scr O_S)$. Equivalently one may
consider the pluridegrees
$$
d_i = (K+H)^i\cdot H^{3-i} = c_{2+i} (\scr J_X(5)),\qquad i=0,\dots, 3,
$$
introduced in \cite{BBS}.
By Zak's theorem \cite{Za} $X$ is linearly normal. Moreover
$h^1(X,\scr O_X) = 0$ by Barth-Larsen-Lefschetz \cite{BL}. In particular $S$
is linearly normal and regular. Clearly $X$ is projectively
Cohen-Macaulay iff $S$ has this property. So by studying $S$ we
obtain from Theorem 4.1 and Theorem 4.2:
\proclaim
{Proposition 5.1}
Let $X$ be contained in a cubic hypersurface. Then $X$ is
projectively Cohen-Macaulay. In particular $X$ is of general type if
$d\geq 13$.\qquad\qed\endproclaim
\noindent To work out a finite list of admissible invariants for a given degree
one may again start with Halphen's upper bound for $\pi$. Further
tools are a congruence obtained from Riemann-Roch \cite{BSS2, 0.11},
the inequalities deduced from the semipositivity of
$\scr N_{X/{\openP^5}}(-1)$ \cite {BOSS1, Proposition 2.2} and adjunction
theory. In the context of section 3 we recall some classification
results. $X_1,\dots,X_{30}$ will denote the 3-folds listed in table 7.3.
The first result follows from Theorem 4.4 and Theorem 4.5 (compare also
\cite{BOSS2}).
\proclaim
{Proposition 5.2}
\roster
\item"{(i)\ }" If $X$ is a scroll over a smooth curve, then $X=X_1$
is a Segre cubic scroll.\qquad
\item"{(ii)\ }" If $X$ is a Fano 3-fold of index 2, then $X=X_2$ is
a complete intersection of two quadric hypersurfaces.
\item"{(iii)\ }" If $X$ is a quadric bundle over a smooth curve, then
$X=X_3$ is a Castelnuovo 3-fold.\qquad\qed
\endroster\endproclaim
\proclaim
{Theorem 5.3} \cite {Ott}.
If $X$ is a scroll over a smooth surface, then $X$ is one of the
following:
\roster
\item"{(i)\ }" a Segre scroll $X=X_1$,
\item"{(ii)\ }" a Bordiga scroll $X=X_4$,
\item"{(iii)\ }" a Palatini scroll $X=X_6$,
\item"{(iv)\ }" a scroll $X=X_{11}$ over a K3 surface.\qquad\qed
\endroster\endproclaim
\noindent From now on we suppose that $X$ is none of the exceptional 3-folds
above. Then the adjunction map $\Phi$ is defined and the connected
morphism $r$ of its Stein factorization contracts the linear
$\openP^2 \subset X$ with normal bundle $\scr O_{\openP^{2}}(-1)$ to points.
\proclaim
{Proposition 5.4} \cite{BSS2}
$r$ is an isomorphism unless $X=X_7$.\qquad\qed
\endproclaim
\noindent From now on we suppose that $X\neq X_7$. Then $X$ coincides
with its first reduction.\newline
The next step in adjunction theory is to study $K+H$. This is big and
nef unless $X$ is one of the special varieties listed in \cite{So5}.
In our case these are classified:
\proclaim
{Theorem 5.5} \cite{BOSS2}.
$K+H$ is big and nef unless
\roster
\item"{(i)\ }" $(X,H)$ is a Fano 3-fold of index 1. Then $X=X_5$ is
a complete intersection of type $(2,3)$.
\item"{(ii)\ }" $(X,H)$ is a Del Pezzo fibration over a smooth curve.
Then $X=X_8$ or $X=X_9$.
\item"{(iii)\ }" $(X,H)$ is a conic bundle over a surface. Then
$X=X_{12}$ or $X=X_{20}$.\qquad\qed
\endroster\endproclaim\noindent
{}From now on we suppose that $K+H$ is big and nef. Then $S$ is of
general type and minimal. Therefore $X$ is called
to be of  {\it log-general type} \cite{BSS1}.
In this case further numerical information is provided by the
generalized Hodge index theorem \cite{BBS, Lemma 1.1} and the
parity relations \cite{BBS, Lemma 1.4}.
%
%
\newline
{}From the Kawamata-Shokurov base point free theorem (see \cite{KMM,
\S 3}) we know that for some $m>0$ the linear system $\vert m(K+H)\vert$
gives rise to a morphism, say $\Psi :X\rightarrow X''$. For $m$ large enough
we can assume that $\Psi$ has connected fibers and normal image. We write
$L'' = \Psi_*(H)$, $K'' = K_{X''}$
and $H'' = K'' + L''$. Then $L''$ and
$H''$ are ample and
$$
\Psi^* (H'') \equiv K + H
$$
(cf. \cite{BFS, (0.2.6)}). $(X'', L'')$ is called the
{\it {second reduction}} of $(X, H)$  \cite{So4},\cite{BFS}.
\proclaim
{Proposition 5.6} \cite{BSS3, Corollary 1.3}
If $d\neq 10$ and $d\neq 13$, then $X''$ is smooth and $\Psi$ is an
isomorphism outside a disjoint union $\scr C$ of smooth curves. Let
$C$ be an irreducible component of $\scr C$ and let $D: =
\Psi^{-1} (C)$. Then the restriction $\Psi_{D}$ of $\Psi$ to
$D$ is a $\openP^1$-bundle $\Psi_{D} : D\rightarrow C$ and $\scr
N_{D\mid F}^X \cong \scr O_{\openP^1}(-1)$ for any fiber $F$ of
$\Psi_{D}$. In fact, $\Psi$ is simply the blowing up along
$\scr C$.\qquad\qed
\endproclaim
\noindent {\bf Remark 5.7.}\quad i) If $d=10$, then there is
exactly one case where $X''$ is not smooth. Namely, for $X=X_{16}$
the second reduction morphism $\Psi$, which is defined by
$\vert K+H\vert$, contracts the quadric surface $K$ to a singular
point $p$. Moreover, $\Psi(X)\subset\openP^6$ is a complete
intersection of type $(2,2,3)$, while $X$ is the
projection from $p$ of $\Psi(X)$ (see also section 7).\par
\noindent ii) From \cite{BSS3, (0.5.1) and (1.1)} and
\cite{Ed, (3.1.3)} it follows that in case
$d=13$ the second reduction is singular iff there exist on $X$
divisors $D\cong\openP^2$, with $\scr N_{D}^X \cong \scr O_{\openP^2}(-2)$,
which are contracted to points. We are not aware of any such example.
\qquad\qed
\medskip\noindent
{\bf Example 5.8.}\quad Let $X\subset\openP^5$ be a smooth 3-fold
with $d=11$ and $\pi = 14$. Then $\chi (\scr O_S) = 8$ and $\chi
(\scr O_X) = 0$\ \ (compare \cite{BSS2}). Every smooth surface in
$\openP^4$ with the same invariants as $S$ is linked $(4,4)$ to a
Castelnuovo surface \cite {Po, Prop. 3.70}. In particular $S$ and
hence $X$ are projectively Cohen-Macaulay with syzygies of type
$$
0\rightarrow 2 \scr O(-5) \oplus \scr O(-6)
\overset\varphi\to\rightarrow 4 \scr O(-4) \rightarrow \scr J_X
\rightarrow 0\ .
$$
Consequently $X=X_{18}$ is linked $(4,4)$ to a Castelnuovo 3-fold.
Conversely this shows the existence of 3-folds of type $X_{18}$
\cite{BSS2}.\par
\noindent What kind of 3-fold is $X$ ?\par
\noindent From the invariants we compute the Kodaira dimension
$\kappa (X) = 0$. In order to show that $X$ is a blown up Calabi-Yau
3-fold we study $|K+H|$. By dualizing $\varphi$ we obtain the resolution
$$
0\leftarrow \omega_X(1) \leftarrow \scr O(1) \oplus 2 \scr O \overset
{^t\varphi}\to\leftarrow 4 \scr O(-1) \leftarrow \scr O(-5)
\leftarrow 0\ .
$$
So $|K+H|$ is base point free, $\dim \vert K+H\vert = 7$, and we
have a well-defined map
$
\Psi_{\mid K+H\mid} : X\rightarrow X''\ ,
$
where $X''$ is a 3-fold in $\openP^7$.
Moreover $h^0(\scr O_S(K_S-H_S))=h^2(\scr O_S(1))=1$ by
Riemann-Roch and Severi's theorem, thus $S$ is minimal and there exists a
rigid curve  $D\in\vert K_S-H_S\vert$ with $H_S\cdot D=4$, $p_a(D)=0$.
In particular, $|K_S|=|D+H_S|$ defines an embedding outside the support of $D$
and maps the divisor $D$ onto a line $L$ in
$\openP^6=\openP(H^0(\scr O_S(K_S)))$. It follows that $\Psi=
\Psi_{\vert K+H\vert} :X\to X''\subset\openP^7$ coincides with the
second reduction morphism. Moreover, by Proposition 5.6, $X''$
is smooth, $K$ is a smooth rational scroll $\openP^1\times\openP^1
\overset {(1,2)}\to\hookrightarrow\openP^5$, which is contracted
by $\Psi$ to the line $L\subset X''$, while $X$ is exactly the blow up of
$X''$ along this line. Riemann-Roch gives $\chi(\scr O_X(2H+2K))=32$,
hence $h^0(\scr I_{X''}(2))\ge h^0(\scr O_{\openP^7}(2))-
h^0(\scr O_X(2H+2K))=4$. In other words, $X''$ lies on 4
linearly independent hyperquadrics. In fact, as one can check,
$X''\subset\openP^7$ is the complete intersection $\Sigma_{(2,2,2,2)}$
of 4 hyperquadrics. Conversely, let $L$
 be a line in $\openP^7$ and
$\Sigma_{(2,2,2,2)}\subset\openP^7$ a smooth complete intersection
of 4 hyperquadrics containing $L$. Then a general projection $X
={\text {proj}}_L\Sigma_{(2,2,2,2)}\subset\openP^5$ is a 3-fold
of type $X_{18}$.\quad\qed
\medskip\noindent
{\bf Remark 5.9.}\quad Similarly, \cite{Po, Proposition 3.59} yields an
easy proof for the uniqueness of the examples of smooth 3-folds
with $d=11$ and $\pi=13$ constructed in \cite{BSS2}. The uniqueness
for the other two families with $d=11$  in \cite{BSS2}
is clear from \cite{GP}.\quad\qed\par\medskip
\noindent The construction via syzygies of all smooth 3-folds
$X\subset\openP^5$
known so far is straightforward. Nevertheless, it is sometimes quite subtle
to determine the structure of $X$. We will give examples of this kind
in the next section.
\head{6. Examples: Two families of birational Calabi-Yau 3-folds in $\openP^5$}
\endhead
\noindent In this section we will construct and study a family of
smooth 3-folds $X\subset \openP^5$ with $d=17$, $\pi =32$, $\chi (\scr
O_X) = 0$ and $\chi (\scr O_S) =24$. We will also describe a family
of smooth 3-folds $X'\subset\openP^5$ obtained via linkage
$X'\underset {(5,6)}\to\sim X$.
\medskip
\noindent Let us first explain how to construct $X$ via syzygies. In
view of Riemann-Roch the following is a plausible Beilinson cohomology
table for $\scr J_X (5)$:
%
$$
\vbox{\offinterlineskip
\halign{
\vrule height10.5pt depth 5.5pt#
&\hbox to 25pt{\hfil$#$\hfil}
&\vrule#
&\hbox to 25pt{\hfil$#$\hfil}
&\vrule#
&\hbox to 25pt{\hfil$#$\hfil}
&\vrule#
&\hbox to 25pt{\hfil$#$\hfil}
&\vrule#
&\hbox to 25pt{\hfil$#$\hfil}
&\vrule#
&\hbox to 25pt{\hfil$#$\hfil}
&\vrule#
&\hbox to 30pt{\hfil}#
&#
\cr
\omit\hbox to 0pt{\hskip-3pt\hbox{$\bigg\uparrow$}\raise5pt\hbox{$i$}\hss}\cr
\multispan{13}\hrulefill\cr
& && && && && && &&\cr
\multispan{13}\hrulefill\cr
&1&& && && && && &\cr
\multispan{13}\hrulefill\cr
& &&4&&2&& && && &\cr
\multispan{13}\hrulefill\cr
& && && && && && &\cr
\multispan{13}\hrulefill\cr
& && && && && && &\cr
\multispan{13}\hrulefill\cr
& && && && && &&2&\cr
\multispan{14}\hrulefill
&\vbox to0pt{\vss\vskip5.5pt\hbox{$\!-\negthickspace\negmedspace@>> m
>$}\vss}\cr
}}
$$
\noindent Suppose that a 3-fold $X$ with these data exists. Then
Beilinson's theorem yields an exact sequence
$$
0\rightarrow \scr F = \scr O(-1)\oplus 4\Omega^4
(4)\overset\varphi\to\rightarrow \scr G = 2\Omega^3(3) \oplus 2\scr O
\rightarrow \scr J_X (5)\rightarrow 0.\tag 6.1
$$
Conversely, as one can check, the minors of a generic
$\varphi\in\Hom (\scr F, \scr G)$ vanish along a smooth 3-fold $X$ as
desired. By construction, $\scr J_X$ has syzygies of type
$$
0\leftarrow \scr J_{X}\leftarrow 2\scr O(-5)\oplus 5\scr O(-6)
\leftarrow 8\scr O (-7) \leftarrow 2\scr O (-8)\leftarrow 0\ .
\tag 6.2
$$
%
%
%
%
%
%
%
What kind of 3-fold is $X$ ? From the syzygies we see that $X$ can be
linked $(5,5)$ to a 3-fold $Z$ of degree $8$. It is not too hard to
identify the scheme $Z$.
\medskip
\noindent Starting conversely with $Z$ we will reconstruct $X$ and
study its geometry. $Z$ can be described as follows: Let
$Y=\openP^1\times\openP^2 \overset {(1,1)}\to\hookrightarrow
\openP^5$ be a Segre cubic scroll and let $L_1,\dots ,L_5$ be five
general lines in $\openP^2$. Then for $i=1,\dots ,5$ the quadric
$Q_i = \openP^1\times L_i \overset {(1,1)}\to\hookrightarrow
\openP^5$ is contained in $Y$ and spans a linear subspace
$\Pi_i\subset\openP^5$ of dimension 3. Clearly, $\Pi_i\cap Y = Q_i$\
($Y$ is cut out by quadrics) and $\Pi_i\cap\Pi_j = \openP^1\times
\{p_{ij}\}$, where $\{p_{ij}\} = L_i\cap L_j$ for $i<j$. Hence the scheme
$$
Z := Y \cup \Pi_1\cup\dots\cup \Pi_5
$$
is locally Cohen-Macaulay, and moreover a local complete intersection
outside the lines $L_{ij} := \openP^1\times \{p_{ij}\}$. Write $Z_k =
Y\cup \bigcup\limits_{i=1}^k \Pi_i$, $k = 0,\dots, 5$. Then $Z_0
= Y$ and $Z_5 = Z$. From the exact sequences
$$
0\rightarrow \scr J_{Z_{k-1}}(m-1)\rightarrow \scr
J_{Z_k}(m)\rightarrow \scr J_{Z_{k-1}\cap \openP^4,\openP^4}(m-1)\rightarrow
0,\tag 6.3
$$
where $\openP^4\subset\openP^5$ is a general hyperplane through $\Pi_k$, we
deduce that $h^0 \scr J_Z (3) = 0$, $h^0\scr J_Z (4) = 1$ and $h^0 \scr J_Z (5)
= 26$,
and that $\scr J_Z (5)$ is globally generated.
\proclaim
{Proposition 6.4}
Let $X$ be linked to $Z$ in the complete intersection of two general
quintic hypersurfaces containing $Z$. Then $X$ is smooth, it contains
the lines $L_{ij}$ and $\scr J_X(5)$ has a resolution of type {\rm (6.1)}.
\endproclaim
\noindent{\bf Proof.}\quad
By a variant of Theorem 2.1 (compare \cite {PS}) $X$ is smooth
outside the lines $L_{ij}$. By using the exact sequences (6.3) we
see that the general quintic hypersurface through $Z$ contains the
first infinitesimal neighborhood of $L_{ij}$, which is a multiplicity
5 structure on such a line. Higher infinitesimal neighborhoods are
not contained in the general quintic hypersurface through $Z$. Moreover the
tangent cone to $Z$ at a point $p\in L_{ij}$ is $\Pi_i \cup \Pi_j\cup \T_p Y$,
and $\Pi_i \cap \T_p Y = \T_p Q_i$. Now a local computation
shows that indeed $X$ is smooth along and contains the
lines $L_{ij}$. That $\scr J_X(5)$ has a Beilinson cohomology table
as above follows via liaison from the exact sequences (6.3).\qquad\qed
\medskip\noindent{\bf Remark 6.5.}\quad
(i)\ By dualizing (6.1) we find that $\omega_X(1)$ has a
presentation of type
$$
0\leftarrow \omega_X(1) \leftarrow \scr O(1) \oplus 18 \scr O
\leftarrow 50 \scr O(-1) \leftarrow \cdots \ .
$$
Thus $|K+H|$ is base point free and $\dim |K+H|=24$.\par\smallskip\noindent
(ii)\ From the double point formulae we compute
$$
H^2\cdot K = 28,\qquad H\cdot K^2 = 18 \text{\quad and\quad} K^3 = -52\ .
$$
In particular $(K+H)\cdot K^2=-34$, hence $\kappa (X)\leq 1$. In fact, as
we will see later, $X$ is a birational Calabi-Yau 3-fold.\qquad\qed
\medskip\noindent
We use in the sequel the above liaison to describe the geometry of $X$:
\proclaim
{Lemma 6.6}
Each linear subspace $\Pi_i$ intersects $X$ along a smooth sextic surface
$S_i$. A
general element in the residual pencil $\vert H-S_i\vert$ is a smooth
blown-up $K3$ surface of degree 11, sectional genus 12, which is embedded in
its
corresponding $\openP^4$ by a linear system of type
$$
\vert H_{min} - 2 E_1 - \sum\limits_{i=2}^{10} E_i\vert\ ,\text{\quad
with\quad} H_{min}^2= 24\ .
$$
\endproclaim
\noindent
{\bf Proof.}\quad It follows from the standard liaison exact sequences
that $X$ meets $\Pi_i$ along a divisor in the class
$4 H_{\Pi_i} - K_{\Pi_i} - Q_i$, hence along a sextic surface $S_i$,
which is smooth for general choices in the liaison. For the
second statement in the lemma, we observe that a general element in
$\vert H-S_i\vert$ is linked $(4,4)$ inside the hyperplane $H$ to
the configuration of planes
$
P_i\cup \bigcup\limits_{j\neq i} (H\cap \Pi_j)\ ,
$
where $P_i$ is the plane residual to $Q_i$ in the intersection $H\cap
Y$. Therefore the lemma follows from the following:
\proclaim
{Proposition 6.7} \cite{Po}
Let $T$ be a configuration $T = P\cup P_1\cup P_2\cup P_3\cup P_4 \subset
\openP^4$,
where $P$ is a plane, while $P_i$, $i=1,\dots,4$, is a plane meeting $P$ along
a line,
such that no three of the lines have common intersection points.
Then $T$ can be linked in the complete intersection of two general quartic
hypersurfaces to a smooth, non-minimal $K3$ surface $S\subset\openP^4$ with
$d=11$ and $\pi =12$, embedded by a linear system
$$
H_S \equiv H_{min} - 2 E_1 - \sum\limits_{i=2}^{10}E_i\ ,\quad H_{min}^2
= 24\ .
$$
Moreover, $P$ meets $S$ along the exceptional conic $E_1$ and an extra
scheme of length 6, while each intersection $P_i\cap S$ is a plane quintic
curve. Residual to it there is a base point free pencil of elliptic space
curves of degree 6.\qquad\qed
\endproclaim
\proclaim
{Lemma 6.8}
$\vert H-S_i\vert$ is a base point free pencil, $i= 1,\dots,5$.
\endproclaim
\noindent{\bf Proof.}\quad
Let $H_i$ denote a general hyperplane containing $\Pi_i$, and let
$K_i$ be the surface residual to $S_i$ in $H_i\cap X$. Then
$S_i\cap K_i\subset\Pi_i\cap K_i$, and in fact equality holds
since $\deg S_i\cap K_i=32-\pi(S_i)-\pi(K_i) +1 = 11$. Thus if
$C_i=S_i\cap K_i$, then $C_i^2=2p_a(C_i)-2-K_{S_i}\cdot C_i= 0$,
where the intersection numbers are computed on $S_i$,
and the lemma follows.\qquad\qed
\medskip\noindent
As a corollary, we deduce that $X$ is an elliptic 3-fold, namely
\proclaim
{Corollary 6.9}
For all $i\neq j$, the linear system $\vert H-S_i\vert \boxtimes
\vert H-S_j\vert$ induces an elliptic fibration
$$
\varphi_{\vert H-S_i\vert \boxtimes \vert H-S_j\vert}
: X \rightarrow \openP^1 \times \openP^1,
$$
with elliptic space curves of degree 6 as fibres.
\endproclaim
\noindent{\bf Proof.}\quad
Fix a general point in $\openP^1\times\openP^1$, i.e., two
general hyperplanes, $H_i$ containing $\Pi_i$ and $H_j$ containing $\Pi_j$,
and denote as above by $K_i$ and $K_j$ the residual surfaces to $S_i$ and
$S_j$ respectively. By Proposition 6.7, $H_i\cap \Pi_j\cap K_i$ and
$H_j\cap \Pi_i\cap K_j$ are plane quintic curves, hence $K_i\cap K_j$
is an elliptic space curve of degree 6, namely the residual to
$H_i\cap\Pi_j\cap K_i$ in $H_j\cap K_i$, or equivalently the residual to
$H_j\cap\Pi_i\cap K_j$ in $H_i\cap K_j$.\qquad\qed
\medskip\noindent
By liaison we deduce that $X$ meets the Segre scroll $Y$ along a
surface $T_2$ in the class $4H_Y- K_Y- \sum\limits_{i=1}^5 Q_i$,
thus along a conic bundle of degree 10 and sectional genus 6.
Moreover, the standard liaison exact sequences yield on $X$ the
linear equivalence
$$
4 H - K\equiv T_2 + \sum\limits_{i=1}^5 S_i.\tag {6.10}
$$
We study in the sequel the structure of the map defined by the composition
of the cartesian product of the 5 pencils $\vert H - S_i\vert$, $i =1,\dots,5$,
with the Segre embedding to $\openP^{31}$:
$$
\Upsilon = \Upsilon_{\vert 5H-\sum\limits_{i=1}^5 S_i\vert}: X
\rightarrow \openP^1 \times\openP^1 \times\openP^1 \times\openP^1
\times\openP^1 \hookrightarrow \openP^{31}\ .
\tag {6.11}
$$
\proclaim
{Lemma 6.12}
The canonical divisor $K$ of $X$ has two components $\T_1$ and $\T_2$.
$\T_1$ is a scroll of degree 18 and sectional genus 10, while
$\T_2$ is the above conic bundle of degree 10 and sectional genus 6.
\endproclaim
\noindent {\bf Proof.}\quad
Let, as in the proof of Lemma 6.8, $K_i$ be a general element in
the pencil $|H - S_i|$. We recall that $(H - S_i)^2 = 0$,
thus $K\mid_{K_i}\equiv K + (H-S_i)\mid_{H-S_i}\equiv K_{H-S_i}$.
In other words, $K$ meets a $K3$ surface $K_i$ along its
canonical divisor, namely, by Proposition 6.7, along 9 exceptional lines and
one
exceptional conic in the plane residual to $Q_i$ in $H\cap Y$.
In conclusion, the exceptional conics sweep the conic bundle $\T_2$,
which is thus a component of $K$, while the exceptional lines on
the $K_i'$s are rulings of a scroll $\T_1$ of degree $H^2\cdot K-10=18$.
Since $\omega_{S_i} = \scr O_{S_i}(2)$ and $S_i\cap S_j=L_{ij}$,
(6.10) restricted to $S_i$ yields  $\T_2\cap S_i
\equiv 2 H_{S_i} - \sum\limits_{j\neq i} L_{ij}$. On the other side
from (6.10) again we infer:
$$
(H-S_i)\mid_{S_i} + (2H - \sum\limits_{j\neq i} S_j)\mid_{S_i} + 2
H\mid_{S_i} \equiv \T_2\mid_{S_i} + H\mid_{S_i} + K\mid_{S_i}.
$$
Thus, since $(H - S_i)\mid_{S_i} \equiv H_{K_i}$, it follows that $K$
intersects $S_i$ along a hyperplane section of $X$.
We deduce that $\T_1$ must intersect $S_i$ along a curve of degree 9
and genus 10, which in turn must be a section of this scroll
since it meets the exceptional lines of $K_i$ in one point.
In other words, $\T_1$ is a scroll of degree 18 and genus 10.\qquad\qed
\proclaim
{Lemma 6.13}
\roster
\item"{i)\ }" The linear system $\vert K+H\vert$ defines a birational
morphism $\Psi=\Psi_{\vert K+H\vert} :X\rightarrow \Psi (X)\subset
\openP^{23}$, which contracts the scroll $T_1$ to a curve of degree 27.
Moreover, $X$ is the blowing up of  $\Psi (X)$ along this curve.
\item"{ii)\ }" The morphism $\Upsilon=\Upsilon_{\vert 5H-\sum\limits_{i=1}^5
S_i\vert}:
X\rightarrow \Upsilon (X)\subset \openP^{31}$, induced by $\vert H+K+T_2\vert$,
contracts the conic bundle $T_2$ to a curve and is birational
on its image.
\endroster
\endproclaim
\noindent {\bf Proof.}\quad
Let $K_i$ be a general element of the pencil $\vert H-S_i\vert$. Part
i) follows easily since $\vert K+H\vert$ induces on $K_i$ the
adjunction morphism $\Phi_i=\Phi_{\vert H_{K_i} + K_{K_i}\vert} :K_i
\rightarrow \Phi_i (K_i)\subset \openP^{12}$, which is birational and
blows down only the 9 exceptional lines $K\cap K_i$. A similar
argument works for part ii) since $\vert 5 H - \sum\limits_{i=1}^5
S_i\vert = \vert H+K+T_2\vert$ restricts to $K_i$ as the map onto
the second adjoint surface given by the adjunction process.\qquad\qed
\medskip\noindent
We can show now that $X$ is a non minimal Calabi-Yau 3-fold, namely:
\proclaim
{Proposition 6.14}
The morphism
$$
\Upsilon = \Upsilon_{\vert 5H-\sum\limits_{i=1}^5 S_i\vert}: X
\rightarrow \openP^1 \times\openP^1 \times\openP^1 \times\openP^1
\times\openP^1 \hookrightarrow \openP^{31}
$$
is birational on its image and contracts only the canonical divisor of $X$
to a curve. Moreover, the image $\Upsilon (X)$ is a smooth complete
intersection of type $(1,1,1,1,1)^2$ in ${\underset{i=1}\to{\overset 5
\to\times}}\openP^1$, hence a minimal Calabi-Yau 3-fold in $\openP^{29}$.
\endproclaim
\noindent {\bf Proof.}\quad The smoothness of $\Upsilon (X)$ follows
from the fact that the iterated adjunction morphisms for $K_i$
blow down only the (-1)-lines and (-1)-conics onto the minimal model of $K_i$.
To see further that $\Upsilon (X)$ is a complete intersection of the type
claimed
we need to compute some intersection numbers. By Lemma 6.8, $(H-S_i)^2 =0$,
thus $(H-S_i)^2\cdot H=0$ and $(H-S_i)^2\cdot S_j = 0$,
which yields $H\cdot S_i^2 = -5$, $S_i^2\cdot S_j = -4$, for $i\neq j$,
and $S_i^3 = -16$. Moreover $S_i\cdot S_j\cdot S_k = 0$, for $i\neq j\neq
k$, $i\neq k$, since $\Pi_i\cap \Pi_j\cap \Pi_k\subset L_{ij}\cap L_{ik}\cap
L_{jk} = \emptyset$, and so we deduce that
$\deg \Upsilon (X) = (5 H-\sum\limits_{i=1}^5 S_i)^3= 120$. On the other
side, $\deg {\underset{i=1}\to{\overset 5 \to\times}}\openP^1=5!=120$
in $\openP^{31}$, while $\Upsilon (X)$ spans only a $\openP^{29}$
since $h^0(\scr O_X(H+K+T_2))\leq h^0 (\scr O_X(H+2K))=
\chi (\scr O_X(H+2K))=30$. The proposition follows.\qquad\qed
\medskip\noindent
\proclaim
{Proposition 6.15}
Let $V^5$ and $V^6$ be general hypersurfaces of degrees 5 and 6 resp.
containing $X$. Then $X$ can be linked in the complete
intersection of $V^5$ and $V^6$ to a smooth 3-fold
$X'\subset\openP^5$. $X'$ has invariants $d' = 13$,
$\pi' = 18$, $\chi (\scr O_{X'}) = 0$, $\chi (\scr O_{S'}) = 10$ and
$p_g(X') = h^0 (\scr J_X(5)) - 1 = 1$. Hence $(H')^2\cdot K' = 8$,
$H'\cdot (K')^2 = -2$ and $(K')^3 = -4$ by the double point formulae.
\endproclaim
\noindent {\bf Proof.}\quad Smoothness follows from a Bertini argument since,
on
$V^5$, $X$ is cut out by sextic hypersurfaces (compare 2.1 and \cite {PS}). The
numerical information follows from the standard liaison exact
sequences.\qquad\qed
\proclaim
{Corollary 6.16}
$X'$ is the degeneracy locus of a morphism
$$
0\rightarrow \scr O (-1) \oplus 2\Omega^2 (2) \rightarrow 4\Omega^1
(1) \oplus 2 \scr O \rightarrow \scr J_{X'}(5) \rightarrow 0\ .
$$
Hence $\scr J_{X'}$ has syzygies of type
$$
0\leftarrow \scr J_{X'}(5) \leftarrow 2\scr O\oplus 19\scr O(-1)
\leftarrow 50\scr O (-2) \leftarrow 48\scr O (-3)
\leftarrow 22\scr O (-4) \leftarrow 4\scr O (-5) \leftarrow 0\ .
$$
\endproclaim
\noindent {\bf Proof.}\quad This follows from $(6.1)$ via
liaison or by applying Beilinson's theorem.\qquad\qed
\medskip\noindent
What type of 3-fold is $X'$?
\proclaim
{Proposition 6.17}
$K'$ is a smooth scroll of degree 8, sectional genus 3 over a
plane quartic curve. Moreover, the Segre scroll
$Y =\openP^1\times\openP^2$ meets $X'$ along the scroll $K'$ and a
curve of degree 9, arithmetic genus 4 on the scroll $T_2$.
\endproclaim
\noindent {\bf Proof.}\quad From general liaison arguments it follows
that $Z$ intersects $X'$ along the canonical divisor of $X'$ plus may be
something of bigger codimension. On the other side, $V^6\cap\Pi_i=S_i$ since
$\Pi_i\cap X=S_i$ for all $i$. We deduce that the 2-dimensional part of
the scheme theoretical intersection $Y\cap X'$ is exactly $K'$.
Now Pic $(Y)$ is generated by the classes $P = [\{\text {point}\}\times
\openP^2]$ and
$Q =[\openP^1 \times \openP^1]$, and $P^2 = 0$, $Q^3 = 0$, $Q^2\cdot P = 1$.
Then the scroll $T_2$ is of class $4 H_Y - K_Y - \sum\limits_{i=1}^5 Q_i\equiv
6P + 2Q$. But $K'$ is residual to $T_2$ in $Y\cap V^6$.
Moreover, $T_2$ is cut out on $Y$ outside the $\Pi_i$'s by the sextic
hypersurfaces through $X$. It follows that $K'$ is smooth for a
general choice of the liaison, and that $K'$ is of class
$6 H_Y- T_2 \equiv 4Q$. In particular, $K'$ is a scroll over a
plane quartic curve and has the claimed invariants. Outside $K'$,\
$X'$ can meet the scroll $Y$ only inside $T_2 \subset Y\cap X$.
The proposition follows now since $X'\cap T_2 \equiv (5 H_X - K_X)\cdot T_2$
is a curve of degree 41, arithmetic genus 80, with a component of degree 32 on
the
scroll $K'$.\qquad\qed
\proclaim
{Proposition 6.18}
\roster
\item"{i)\ }" $\vert H'+K'\vert$ is base point free and big.
\item"{ii)\ }" $\Psi'=\Psi_{\vert H'+K'\vert}:X'\rightarrow \openP^9$ is
birational on its image $M = \Psi' (X')$, which is
a smooth Calabi-Yau 3-fold, with $\deg M = 27$, $\pi (M) = 28$ and
$c_3 (M) = -64$.
\item"{iii)\ }" $\Psi'$ contracts the scroll
$K'$ to a curve of degree 6 and genus 3, and is an
isomorphism outside this scroll. Moreover, $X'$ is the blow up of
$M$ along this curve.
\endroster
\endproclaim
\noindent {\bf Proof.}\quad
i) From the syzygies we see that $\omega_{X'}$ is a quotient
of $\scr O\oplus 4\scr O(-1)$, thus $\vert H'+K'\vert$ is base point
free and big since $(H'+K')^3 =27$. Moreover $\dim\vert H'+K'\vert=9$
\par\smallskip\noindent
ii) From the liaison exact sequence
$$
0\rightarrow \scr J_{V^5\cap V^6} (6)\rightarrow \scr J_X (6) \rightarrow
\omega_{X'} (1) \rightarrow 0
$$
we deduce that the map $\Psi': X' \rightarrow
\openP^9$ is in fact the composition of the restriction to $X'$ of
the rational morphism $\Xi : \openP^5 \dashrightarrow \openP^{16}=
\openP (H^0(\scr J_X(6))$ given by the sextic hypersurfaces through $X$,
with a projection from $\openP^{16} \dashrightarrow \openP^9$ along
$\openP^6 =\openP (H^0(\scr J_{V^5\cap V^6}(6))$. Thus in order to show that
$\Psi'$ is birational on its image it is enough to check that
$\Xi$ is birational on its image and that the projection $\openP^{16}
\dashrightarrow \openP^9$ is generic. But one checks easily that 5
general sextic hypersurfaces through $X$ meet in exactly one point
outside $X$. In particular, it follows that $\Psi':
X'\rightarrow M$ coincides with the second reduction map of $X'$.
\par\smallskip\noindent
iii) Since $(H' + K')\cdot K' = 0$, $\Psi'$
contracts the scroll $K'$. Its image is isomorphic to the plane
quartic curve, which is the base of the scroll $K'$. From Remark 5.7 it follows
that $M$ is smooth and $\Psi'$ is an isomorphism
outside the scroll $K'$, unless there are divisors
$D\cong\openP^2\subset X'$ with ${\scr N}_{D}^{X'}\cong\scr O_{\openP^2}(-2)$,
which are contracted to singular points on $M$.
Assume that such a divisor $D$ exists. Then a general hyperplane
section $S'$ of $X'$ contains a $(-2)$-line $L$. But on $S'$,
$h^0(\scr O_{S'} (K_{S'} - H_{S'})) = h^3 (\scr O_{X'}) = 1$, so if
$D_{S'}\in\vert K_{S'} - H_{S'}\vert$, then $K_{S'}\cdot L = 0 =
1 + D_{S'}\cdot L$, and thus $L$ must be a component of $D_{S'}$.
On the other side  $D_{S'}$ is an irreducible hyperplane section of the
smooth scroll $K_{X'}$, and therefore contains no such line as a component.
\qquad\qed

\vfil\eject

\head{7. Overview}\endhead
\noindent
In this section we collect some information on the families of smooth
non general type surfaces in $\openP^4$ and 3-folds in $\openP^5$ known to us.
%
$$\gather
\hbox{\text Table 7.1.\quad Known families of smooth
non general type surfaces in $\openP^4$ } \\
\vbox{\tabskip=0pt \offinterlineskip
\def\tablerule{\noalign{\hrule}}
\halign to391pt{\strut#& \vrule#\tabskip= 1em plus2em &\hfil#&
\vrule#& \hfil# \hfil&
\vrule#& \hfil# \hfil&
\vrule#& \hfil# \hfil&
\vrule#& \hfil# \hfil&
\vrule#& \hfil# \hfil&
\vrule#& \hfil# \hfil&
\vrule#& \hfil# &\vrule#
\tabskip=0pt\cr\tablerule\tablerule
&&&&\multispan{13}\hfil Enriques-Kodaira Classification\hfil&
\cr\tablerule\tablerule
&&\omit\hidewidth degree\hidewidth&&\omit\hidewidth rational\hidewidth
&&\omit\hidewidth ruled\hidewidth&&\omit\hidewidth Enriques\hidewidth
&&\omit\hidewidth {$K3$}\hidewidth&&{abelian}
&&\omit\hidewidth
{bielliptic}\hidewidth&&\omit\hidewidth{elliptic}\hidewidth&\cr
&&\omit&&\omit&&\omit\hidewidth irrat.\hidewidth&&\omit&&\omit&&\omit&&
\omit&&\omit&\cr\tablerule\tablerule
%
\tablerule\tablerule
&&{$d\le 4$}&&6&&\omit&&\omit&&1&&\omit&&\omit&&\omit&\cr\tablerule\tablerule
&&{$d=5$}&&1&&1&&\omit&&\omit&&\omit&&\omit&&\omit&\cr
&&\omit&&[Ca]&&\omit\hidewidth[Seg]\hidewidth&&\omit&&\omit&&\omit&&
\omit&&\omit&\cr\tablerule\tablerule
&&{$d=6$}&&1&&\omit&&\omit&&1&&\omit&&\omit&&\omit&\cr
&&\omit\hidewidth [Io1],[Ok1]\hidewidth&&\omit\hidewidth[Bo],[Ve]\hidewidth
&&\omit&&\omit&&\omit&&\omit&&\omit&&\omit&\cr
&&\omit&&\omit[Wh]&&\omit&&\omit&&\omit&&
\omit&&\omit&&\omit&\cr\tablerule\tablerule
&&{$d=7$}&&1&&\omit&&\omit&&1&&\omit&&\omit&&1&\cr
&&\omit\hidewidth [Io1],[Ok3]\hidewidth
&&\omit\hidewidth [Io1],[Ok3]\hidewidth
&&\omit&&\omit&&[Ro1]&&\omit&&\omit&&
\omit\hidewidth[Ba]\hidewidth&\cr\tablerule\tablerule
&&{$d=8$}&&2&&\omit&&\omit&&1&&\omit&&\omit&&1&\cr
&&\omit\hidewidth [Io2],[Ok4]\hidewidth&&\omit\hidewidth[Ok4],[Al1]\hidewidth
&&\omit&&\omit&&[Ok4]&&\omit&&\omit&&\omit\hidewidth[Ba]\hidewidth&\cr
\tablerule\tablerule
&&{$d=9$}&&2&&\omit&&1&&1&&\omit&&\omit&&1&\cr
&&[AR]&&\omit\hidewidth[Al1],[Al2]\hidewidth&&\omit&&\omit\hidewidth[Cos],[CV]
\hidewidth&&[Ro1]&&\omit&&\omit&&\omit\hidewidth[AR]\hidewidth&\cr
\tablerule\tablerule
&&\omit\hidewidth{$d=10$}\hidewidth&&3&&\omit&&1&&2&&1&&1&&2&\cr
&&\omit\hidewidth[Ra],[PR]\hidewidth&&\omit\hidewidth[DES],[Ra]\hidewidth&&
\omit&&\omit\hidewidth[DES],[Br]\hidewidth&&\omit\hidewidth[Ra],[Po]\hidewidth
&&\omit\hidewidth[Co],[HM]\hidewidth&&[Ser]&&\omit\hidewidth[Ra]
\hidewidth&\cr
&&\omit&&\omit&&\omit&&\omit&&\omit&&\omit\hidewidth[HL],[Ram]\hidewidth
&&\omit\hidewidth[ADHPR1]\hidewidth&&\omit&\cr\tablerule\tablerule
&&\omit\hidewidth{$d=11$}\hidewidth&&\omit\hidewidth 3+2\hidewidth
&&\omit&&1&&5&&\omit&&\omit&&1&\cr
&&[Po]&&\omit\hidewidth[DES],[Po]\hidewidth&&\omit&&[DES]&&\omit
\hidewidth[DES],[Po]\hidewidth&&\omit&&\omit&&\omit
\hidewidth[Po]\hidewidth&\cr\tablerule\tablerule
&&\omit\hidewidth{$d=12$}\hidewidth&&\omit
&&\omit&&\omit&&1&&\omit&&\omit&&3&\cr
&&\omit&&\omit&&\omit&&\omit&&\omit\hidewidth[DES]\hidewidth&&\omit&&\omit&&
\omit\hidewidth[Po]\hidewidth&\cr\tablerule\tablerule
&&\omit\hidewidth{$d=13$}\hidewidth&&\omit
&&\omit&&2&&1&&\omit&&\omit&&\omit&\cr
&&\omit&&\omit&&\omit&&\omit\hidewidth[DES],[Po]\hidewidth&&
[Po]&&\omit&&\omit&&\omit&\cr\tablerule\tablerule
&&\omit\hidewidth{$d=14$}\hidewidth&&\omit&&\omit
&&\omit&&1&&\omit&&\omit&&\omit&\cr
&&\omit&&\omit&&\omit&&\omit&&[Po]&&\omit&&
\omit&&\omit&\cr\tablerule\tablerule
&&\omit\hidewidth{$d=15$}\hidewidth&&\omit
&&\omit&&\omit&&\omit&&2&&1&&\omit&\cr
&&\omit&&\omit&&\omit&&\omit&&\omit&&\omit\hidewidth[HM],[Po]\hidewidth
&&\omit\hidewidth[ADHPR1]\hidewidth&&\omit&\cr\tablerule\tablerule
\hfil\cr}} \endgather$$
%
\par\medskip\noindent
{\bf Remark 7.2.}\quad (i) The classification of smooth surfaces in
$\openP^4$ is complete up to degree 10, and there is a partial
classification in degree 11. In the first column of Table 7.1
we refer to the papers, where one can find the classification results.
In the other columns we indicate the number of families known and the
corresponding references. The classification up to degree 5 is
classical. More information can be found in \cite{DES, Appendix B},
\cite{Po, Appendix} and \cite{ADHPR2}.\par\noindent
(ii) Two families of rational surfaces of degree 11 are due to
Schreyer (unpublished).\par\noindent
(iii) One of the families of K3 surfaces of degree 11 has been first
constructed by Ranestad (compare \cite{Po, Proposition 3.41}).\quad\qed
\par\medskip
\vfil\eject
\noindent
\def\link#1#2{{\smash{\mathop{\sim}\limits^{\scriptscriptstyle(#1,#2)}}}}
\def\un{\mathop{\cup}}
\def\kreuz{\mathop{\times}}
$$\gather
\hbox{\text Table 7.3.\quad Known families of smooth, non-degenerate,
non general type 3-folds in $\openP^5$ } \\
{\eightpoint
\vbox{\tabskip=0pt\offinterlineskip
\def\tablerule{\noalign{\hrule}}
\halign {\strut#& \vrule#\tabskip=1em plus2em&
\hfil#
&\vrule#&\hfil#\hfil&\vrule#&\hfil#\hfil
&\vrule#&\hfil#\hfil&\vrule#&\hfil#\hfil
&\vrule#&\hfil#\hfil&\vrule#&\hfil#\hfil
&\vrule#&\hfil#\hfil&\vrule#&\hfil#\hfil
&\vrule#&\hfil#\hfil&\vrule#&\hfil#
&\vrule#\tabskip=0pt\cr\tablerule\tablerule
%
%
&&\omit\hidewidth $X$\hidewidth&&
\omit\hidewidth $d$\hidewidth&&
\omit\hidewidth $\pi$\hidewidth&&
\omit\hidewidth $p_g$\hidewidth&&
\omit\hidewidth $\chi({\scr O}_X)$\hidewidth&&
\omit\hidewidth $\chi({\scr O}_S)$\hidewidth&&
\omit\hidewidth $\kappa(X)$\hidewidth&&
\omit\hidewidth liaison\hidewidth&&
\omit\hidewidth type\hidewidth&&
\omit\hidewidth classification\hidewidth&&
\omit\hidewidth ref.\hidewidth&
\cr\tablerule\tablerule
&&\omit&&\omit&&\omit&&\omit&&\omit&&\omit
&&\omit&&\omit&&\omit&&\omit&&\omit&\cr
&&\omit\hidewidth$X_1$\hidewidth&&3&&0&&0&&1&&1&&\omit\hidewidth
$-\infty$\hidewidth
&&$X_1\link 2 2\openP^3$
&&\omit\hidewidth Segre embedd. of
$\openP^1\times\openP^2$\hidewidth&&rational&&\omit&\cr
&&\omit&&\omit&&\omit&&\omit&&\omit&&\omit
&&\omit&&\omit&&\omit&&scroll&&\omit&\cr
\tablerule\tablerule
&&\omit&&\omit&&\omit&&\omit&&\omit&&\omit
&&\omit&&\omit&&\omit&&\omit&&\omit&\cr
&&\omit\hidewidth$X_2$\hidewidth&&4&&1&&0&&1&&1&&\omit\hidewidth
$-\infty$\hidewidth
&&\omit\hidewidth$X_2=\Sigma_{(2,2)}$\hidewidth&&Fano 3-fold of index
2&&rational
&&[Kl]&\cr
&&\omit&&\omit&&\omit&&\omit&&\omit&&\omit
&&\omit&&\omit&&\omit&&\omit&&\omit&\cr
\tablerule\tablerule
&&\omit&&\omit&&\omit&&\omit&&\omit&&\omit&&\omit&&\omit&&Castelnuovo
3-fold;&&\omit&&\omit&\cr
&&\omit\hidewidth$X_3$\hidewidth&&5&&2&&0&&1&&1&&\omit\hidewidth
$-\infty$\hidewidth
&&\omit\hidewidth $X_3\link 2 3\openP^3$\hidewidth
&&\omit\hidewidth quadric fibration over $\openP^1$\hidewidth&&rational
&&\omit\hidewidth [Io1],[Ok2]\hidewidth&\cr
&&\omit&&\omit&&\omit&&\omit&&\omit&&\omit&&\omit&&\omit&&via
$|K+2H|$&&\omit&&\omit&\cr
\tablerule\tablerule
&&\omit&&\omit&&\omit&&\omit&&\omit&&\omit&&\omit&&\omit&&Bordiga
3-fold;&&\omit&&\omit&\cr
&&\omit\hidewidth$X_4$\hidewidth&&6&&3&&0&&1&&1&&\omit\hidewidth
$-\infty$\hidewidth
&&\omit\hidewidth$X_4\link 3 3 X_1$\hidewidth
&&$X_4=\openP({\scr E})$, ${\scr E}$ rk 2 vb on $\openP^2$
&&rational&&\omit\hidewidth [Io1],[Ok2]\hidewidth&\cr
&&\omit&&\omit&&\omit&&\omit&&\omit&&\omit&&\omit&&\omit
&&\omit\hidewidth$c_1=4$, $c_2=10$, via $|K+2H|$\hidewidth
&&scroll&&\omit&\cr
\tablerule\tablerule
&&\omit&&\omit&&\omit&&\omit&&\omit&&\omit&&\omit&&\omit&&\omit
&&\omit\hidewidth unirational\hidewidth&&\omit&\cr
&&\omit\hidewidth$X_5$\hidewidth&&6&&4&&0&&1&&2&&\omit\hidewidth
$-\infty$\hidewidth
&&\omit\hidewidth$X_5=\Sigma_{(2,3)}$\hidewidth&&Fano 3-fold of index 1&&not
&&\omit\hidewidth [Io1],[Ok2]\hidewidth&\cr
&&\omit&&\omit&&\omit&&\omit&&\omit&&\omit&&\omit&&\omit&&\omit
&&\omit\hidewidth rational\hidewidth
&&\omit\hidewidth [En],[Fa1]\hidewidth&\cr
\tablerule\tablerule
&&\omit&&\omit&&\omit&&\omit&&\omit&&\omit
&&\omit&&\omit&&\omit&&\omit&&\omit&\cr
&&\omit\hidewidth$X_6$\hidewidth&&7&&4&&0&&1&&1&&\omit\hidewidth
$-\infty$\hidewidth
&&\omit\hidewidth$X_6$ proj.\hidewidth
&&\omit\hidewidth $X_6=\openP({\scr E})$, ${\scr E}$ rk 2 vb on \hidewidth
&&rational&&\omit\hidewidth[Io1],[Ok3]\hidewidth&\cr
&&\omit&&\omit&&\omit&&\omit&&\omit&&\omit&&\omit
&&\omit\hidewidth Buchsbaum\hidewidth
&&\omit\hidewidth $\openP^2(x_1,...x_6)$; via $|K+2H|$\hidewidth
&&scroll&&[Pa]&\cr
\tablerule\tablerule
&&\omit&&\omit&&\omit&&\omit&&\omit&&\omit&&\omit&&\omit
&&\omit\hidewidth$X_7=\Sigma_{(2,2,2)}(x_0)$, the blowing\hidewidth
&&\omit\hidewidth blown up\hidewidth&& &\cr
&&\omit\hidewidth$X_7$\hidewidth&&7&&5&&0&&1&&2&&\omit\hidewidth
$-\infty$\hidewidth
&&\omit\hidewidth$X_7\link 3 3\Sigma_{(1,2)}$\hidewidth
&&\omit\hidewidth up of a c.i. $\Sigma_{(2,2,2)}\subset\openP^6$;\hidewidth
&&\omit\hidewidth Fano of\hidewidth&&\omit\hidewidth[Io1],[Ok3]\hidewidth&\cr
&&\omit&&\omit&&\omit&&\omit&&\omit&&\omit&&\omit&&\omit
&&via $|K+2H|$&&\omit\hidewidth index 1\hidewidth&&\omit&\cr
\tablerule\tablerule
&&\omit&&\omit&&\omit&&\omit&&\omit&&\omit&&\omit&&\omit
&&\omit\hidewidth Del Pezzo fibration over
$\openP^1$,\hidewidth&&\omit&&\omit&\cr
&&\omit\hidewidth$X_8$\hidewidth&&7&&6&&0&&1&&3&&\omit\hidewidth
$-\infty$\hidewidth
&&\omit\hidewidth$X_8\link 2 4\openP^3$\hidewidth
&&\omit\hidewidth with gen. fibre $\openP^2(x_1,...x_6)$;\hidewidth
&&rational&&\omit\hidewidth [Io1],[Ok3]\hidewidth&\cr
&&\omit&&\omit&&\omit&&\omit&&\omit&&\omit&&\omit&&\omit&&via
$|K+H|$&&\omit&&\omit&\cr
\tablerule\tablerule
&&\omit&&\omit&&\omit&&\omit&&\omit&&\omit&&\omit&&\omit
&&\omit\hidewidth Del Pezzo fibration over
$\openP^1$;\hidewidth&&\omit&&\omit&\cr
&&\omit\hidewidth$X_9$\hidewidth&&8&&7&&0&&1&&3&&\omit\hidewidth
$-\infty$\hidewidth
&&\omit\hidewidth$X_8\link 3 3\openP^3$\hidewidth
&&\omit\hidewidth gen. fibre c.i. $(2,2)$ in $\openP^4$;&&rational&&[Io2]&\cr
&&\omit&&\omit&&\omit&&\omit&&\omit&&\omit&&\omit&&\omit&&via
$|K+H|$&&\omit&&\omit&\cr
\tablerule\tablerule
&&\omit&&\omit&&\omit&&\omit&&\omit&&\omit
&&\omit&&\omit&&\multispan{3}&&\omit&\cr
&&\omit\hidewidth$X_{10}$\hidewidth&&8&&9&&1&&0&&6&&$0$
&&\omit\hidewidth$X_{10}=\Sigma_{(2,4)}$\hidewidth
&&\multispan{3}\hfil minimal Calabi-Yau 3-fold\hfil&&\omit&\cr
&&\omit&&\omit&&\omit&&\omit&&\omit&&\omit&&\omit
&&\omit&&\multispan{3}&&\omit&\cr
\tablerule\tablerule
&&\omit&&\omit&&\omit&&\omit&&\omit&&\omit&&\omit&&\omit
&&\omit\hidewidth $\openP^1$ bundle over a
minimal\hidewidth&&scroll,&&\omit&\cr
&&\omit\hidewidth$X_{11}$\hidewidth&&9&&8&&0&&2&&2
&&\omit\hidewidth $-\infty$\hidewidth
&&\omit\hidewidth$X_{11}$ proj.\hidewidth
&&\omit\hidewidth $K3$ surface $S\subset\openP^8$;\hidewidth&&not&&[Ch3]&\cr
&&\omit&&\omit&&\omit&&\omit&&\omit&&\omit&&\omit
&&\omit\hidewidth Buchsbaum\hidewidth&&via $|K+H|$&&rational&&\omit&\cr
\tablerule\tablerule
&&\omit&&\omit&&\omit&&\omit&&\omit&&\omit&&\omit&&\omit
&&\omit\hidewidth conic bundle over\hidewidth&&\omit&&\omit&\cr
&&\omit\hidewidth$X_{12}$\hidewidth&&9&&9&&0&&1&&4&&\omit\hidewidth
$-\infty$\hidewidth
&&\omit\hidewidth$X_{12}\link 3 4X_3$\hidewidth
&&\omit\hidewidth$\openP^2$, via $|K+H|$\hidewidth
&&rational&&\omit\hidewidth [BSS1]\hidewidth&\cr
&&\omit&&\omit&&\omit&&\omit&&\omit&&\omit&&\omit
&&\omit&&\omit&&\omit&&\omit&\cr
\tablerule\tablerule
&&\omit&&\omit&&\omit&&\omit&&\omit&&\omit
&&\omit&&\omit&&\multispan{3}&&\omit&\cr
&&\omit\hidewidth$X_{13}$\hidewidth&&9&&\omit\hidewidth
10\hidewidth&&1&&0&&6&&0
&&\omit\hidewidth$X_{13}=\Sigma_{(3,3)}$\hidewidth
&&\multispan{3}\hfil minimal Calabi-Yau 3-fold\hfil&&\omit&\cr
&&\omit&&\omit&&\omit&&\omit&&\omit&&\omit&&\omit
&&\omit&&\multispan{3}&&\omit&\cr
\tablerule\tablerule
&&\omit&&\omit&&\omit&&\omit&&\omit&&\omit&&\omit&&\omit
&&\omit\hidewidth minimal $K3$ fibration over\hidewidth&&\omit&&\omit&\cr
&&\omit\hidewidth$X_{14}$\hidewidth&&9&&\omit\hidewidth 12\hidewidth&&2
&&\omit\hidewidth -1\hidewidth&&9&&1
&&\omit\hidewidth$X_{14}\link 2 5\openP^3$\hidewidth
&&\omit\hidewidth$\openP^1$, via $|K|$; $4K+3H>0$.\hidewidth
&&\omit&&\omit\hidewidth [BSS1]\hidewidth&\cr
&&\omit&&\omit&&\omit&&\omit&&\omit&&\omit&&\omit&&\omit
&&\omit\hidewidth log-gen type\hidewidth&&\omit&&\omit&\cr
\tablerule\tablerule
&&\omit&&\omit&&\omit&&\omit&&\omit&&\omit&&\omit&&\omit
&&\omit\hidewidth log-general type;\hidewidth&&\omit&&\omit&\cr
&&\omit\hidewidth$X_{15}$\hidewidth&&\omit\hidewidth 10\hidewidth
&&\omit\hidewidth 11\hidewidth&&0&&1&&5&&\omit\hidewidth$-\infty$\hidewidth
&&\omit\hidewidth$X_{15}\link 4 4 X_4$\hidewidth
&&\omit\hidewidth $|K+H|$ is birational onto $\openP^3$\hidewidth
&&rational&&\omit\hidewidth [BSS1]\hidewidth&\cr
&&\omit&&\omit&&\omit&&\omit&&\omit&&\omit&&\omit
&&\omit&&\omit&&\omit&&\omit&\cr
\tablerule\tablerule
&&\omit&&\omit&&\omit&&\omit&&\omit&&\omit&&\omit&&\omit
&&\omit\hidewidth$X_{16}={\text {proj}}_p\Sigma_{(2,2,3)}$, with\hidewidth
&&\omit\hidewidth birational\hidewidth&&\omit&\cr
&&\omit\hidewidth$X_{16}$\hidewidth&&\omit\hidewidth 10\hidewidth
&&\omit\hidewidth 12\hidewidth&&1&&0&&7&&0
&&\omit\hidewidth$X_{16}\link 3 4\Sigma_{(1,2)}$\hidewidth
&&\omit\hidewidth $\Sigma_{(2,2,3)}\subset\openP^6$ c.i. singular at
$p$;\hidewidth
&&\omit\hidewidth Calabi-Yau;\hidewidth&&\omit\hidewidth [BSS1]\hidewidth&\cr
&&\omit&&\omit&&\omit&&\omit&&\omit&&\omit&&\omit&&\omit
&&\omit\hidewidth via the inverse of $|K+H|$\hidewidth
&&\omit\hidewidth $H^2K=2$\hidewidth&&\omit&\cr
\tablerule\tablerule
&&\omit&&\omit&&\omit&&\omit&&\omit&&\omit&&\omit&&\omit
&&\omit\hidewidth blown up Fano 3-fold;\hidewidth
&&\omit\hidewidth unirational\hidewidth&&\omit&\cr
&&\omit\hidewidth$X_{17}$\hidewidth&&\omit\hidewidth 11\hidewidth
&&\omit\hidewidth 13\hidewidth&&0&&1&&6&&\omit\hidewidth $-\infty$\hidewidth
&&\omit\hidewidth$X_{17}\link 4 5 X_{11}$\hidewidth
&&\omit\hidewidth$|K+H|$ is birational onto a \hidewidth&&not&&[Ch3]&\cr
&&\omit&&\omit&&\omit&&\omit&&\omit&&\omit&&\omit
&&\omit&&hypercubic in $\openP^4$;
&&rational&&\omit\hidewidth [BSS2]\hidewidth&\cr
\tablerule\tablerule
&&\omit&&\omit&&\omit&&\omit&&\omit&&\omit&&\omit&&\omit
&&\omit\hidewidth$X_{18}={\text {proj}}_L\Sigma_{(2,2,2,2)}$, with\hidewidth
&&\omit\hidewidth birational\hidewidth&&\omit&\cr
&&\omit\hidewidth$X_{18}$\hidewidth&&\omit\hidewidth
11\hidewidth&&\omit\hidewidth 14\hidewidth
&&1&&0&&8&&0&&\omit\hidewidth$X_{18}\link 4 4 X_3$\hidewidth
&&\omit\hidewidth $L\subset\Sigma_{(2,2,2,2)}\subset\openP^7$ smooth
c.i.\hidewidth
&&\omit\hidewidth Calabi-Yau;\hidewidth&&\omit\hidewidth [BSS2]\hidewidth&\cr
&&\omit&&\omit&&\omit&&\omit&&\omit&&\omit&&\omit&&\omit
&&\omit\hidewidth and $L$ a line; via $|K+H|$\hidewidth
&&\omit\hidewidth$H^2K=4$\hidewidth&&\omit&\cr
\tablerule\tablerule
&&\omit&&\omit&&\omit&&\omit&&\omit&&\omit&&\omit&&\omit
&&\omit\hidewidth minimal $K3$ fibration over\hidewidth
&&\omit&&\omit&\cr
&&\omit\hidewidth$X_{19}$\hidewidth&&\omit\hidewidth 11\hidewidth
&&\omit\hidewidth 15\hidewidth&&2&&\omit\hidewidth -1\hidewidth
&&\omit\hidewidth 10\hidewidth&&1&&\omit\hidewidth$X_{19}\link 3
4\openP^3$\hidewidth
&&\omit\hidewidth $\openP^1$ via $|K|$; fibres are $(2,3)$\hidewidth
&&\omit&&\omit\hidewidth [BSS2]\hidewidth&\cr
&&\omit&&\omit&&\omit&&\omit&&\omit&&\omit&&\omit&&\omit
&&\omit\hidewidth c.i. in $\openP^4$\hidewidth&&\omit&&\omit&\cr
\tablerule\tablerule
\hfil\cr}}}\endgather
$$
\vfil\eject
\noindent
$$\gather
{\eightpoint
\vbox{\tabskip=0pt \offinterlineskip
\def\tablerule{\noalign{\hrule}}
\halign {\strut#&\vrule#\tabskip= 1em plus2em &
\hfil#
&\vrule#&\hfil#\hfil
&\vrule#&\hfil#\hfil
&\vrule#&\hfil#\hfil
&\vrule#&\hfil#\hfil
&\vrule#&\hfil#\hfil
&\vrule#&\hfil#\hfil
&\vrule#&\hfil#\hfil
&\vrule#&\hfil#\hfil
&\vrule#&\hfil#\hfil
&\vrule#&\hfil#&\vrule#\tabskip=0pt\cr
\tablerule\tablerule
&&\omit\hidewidth $X$\hidewidth&&
$d$&&
\omit\hidewidth $\pi$\hidewidth&&
\omit\hidewidth $p_g$\hidewidth&&
\omit\hidewidth $\chi({\scr O}_X)$\hidewidth&&
\omit\hidewidth $\chi({\scr O}_S)$\hidewidth&&
\omit\hidewidth $\kappa(X)$\hidewidth&&
\omit\hidewidth liaison\hidewidth&&
\omit\hidewidth type\hidewidth&&
\omit\hidewidth classification\hidewidth&&
\omit\hidewidth ref\hidewidth&
\cr\tablerule\tablerule
&& && && && && && && && &&conic bundle over a&&not&& &\cr
&&\omit\hidewidth$X_{20}$\hidewidth
&&\omit\hidewidth 12\hidewidth&&\omit\hidewidth 15\hidewidth
&&0&&2&&6&&\omit\hidewidth $-\infty$\hidewidth
&&\omit\hidewidth$X_{20}\link 5 5 X_{25}$\hidewidth
&&$K3$ quartic surface $S\subset\openP^3$,&&rational
&&\omit\hidewidth [BOSS2]\hidewidth&\cr
&& && && && && && && && &&via $|K+H|$&& && &\cr
\tablerule\tablerule
&& && && && && && && &&
&&$|K+H|$ defines a birational&& && &\cr
&&\omit\hidewidth$X_{21}$\hidewidth&&\omit\hidewidth 12\hidewidth
&&\omit\hidewidth 15\hidewidth&&0&&1&&7&&\omit\hidewidth $-\infty$\hidewidth
&&\omit\hidewidth $X_{21}\link 5 5$\hidewidth
&&\omit\hidewidth map onto a Bordiga $X_4\subset\openP^5$\hidewidth
&&rational&&[Ed]&\cr
&& && && && && && && &&\omit\hidewidth$X_{15}\cup X_1$\hidewidth&& && && &\cr
\tablerule\tablerule
&& && && && && && && &&
&&\omit\hidewidth $|K+H|$ defines a birational\hidewidth
&&\omit\hidewidth birational\hidewidth&& &\cr
&&\omit\hidewidth$X_{22}$\hidewidth
&&\omit\hidewidth 12\hidewidth&&\omit\hidewidth 16\hidewidth&&1&&0&&9&&0
&&\omit\hidewidth $X_{22}\link 5 5 X_{27}$\hidewidth&&
\omit\hidewidth morphism onto $\openP(R_2)\cap
Bl_{\openP^2}\openP^8$;\hidewidth
&&\omit\hidewidth Calabi-Yau;\hidewidth&&[Ch3]&\cr
&& && && && && && && &&\omit
&&\omit\hidewidth see [Ch3] for details\hidewidth
&&\omit\hidewidth$H^2K=6$\hidewidth&& &\cr
\tablerule\tablerule
&& && && && && && && &&
&&\omit\hidewidth $|K|$ defines a $K3$ fibration\hidewidth&& && &\cr
&&\omit\hidewidth$X_{23}$\hidewidth&&\omit\hidewidth 12\hidewidth
&&\omit\hidewidth 17\hidewidth&&2
&&\omit\hidewidth -1\hidewidth&&\omit\hidewidth 11\hidewidth&&1
&&$X_{23}\link 4 5 X_9$
&&\omit\hidewidth over $\openP^1$; the fibres are\hidewidth
&&\omit\hidewidth minimal\hidewidth&&[Ed]&\cr
&& && && && && && && && &&$\Sigma_{(2,2,2)}$ in $\openP^5$&& && &\cr
\tablerule\tablerule
&& && && && && && && &&
&&\omit\hidewidth $|K|$ defines an elliptic fibration\hidewidth&& && &\cr
&&\omit\hidewidth$X_{24}$\hidewidth&&\omit\hidewidth 12\hidewidth
&&\omit\hidewidth 18\hidewidth&&3
&&\omit\hidewidth -2\hidewidth&&\omit\hidewidth 13\hidewidth&&2
&&\omit\hidewidth$X_{24}\link 3 6 X_4$\hidewidth&&
\omit\hidewidth (in plane cubics) over $\openP^2$\hidewidth
&&\omit\hidewidth minimal\hidewidth&&[Ed]&\cr
&& && && && && && && && && &&elliptic&& &\cr
\tablerule\tablerule
&& && && && && && && && &&\omit\hidewidth log-gen type\hidewidth&& && &\cr
&&\omit\hidewidth$X_{25}$\hidewidth&&\omit\hidewidth 13\hidewidth
&&\omit\hidewidth 18\hidewidth&&0&&1&&9&&\omit\hidewidth $-\infty$\hidewidth
&&\omit\hidewidth$X_{25}\link 5 5 X_{20}$\hidewidth
&& && &&\omit\hidewidth [BOSS2]\hidewidth&\cr
&& && && && && && && && && && && &\cr
\tablerule\tablerule
&& && && && && && && &&
&&\omit\hidewidth $|K+H|$ defines a birational\hidewidth
&&\omit\hidewidth blown up\hidewidth&& &\cr
&&\omit\hidewidth$X_{26}$\hidewidth&&\omit\hidewidth 13\hidewidth
&&\omit\hidewidth 18\hidewidth&&1&&0&&\omit\hidewidth 10\hidewidth&&0
&&\omit\hidewidth$X_{26}\link 5 6 X_{29}$\hidewidth
&&\omit\hidewidth map onto a Calabi-Yau 3-fold\hidewidth&&
\omit\hidewidth Calabi-Yau\hidewidth&& &\cr
&& && && && && && && &&
&&\omit\hidewidth $Y\subset\openP^9$ with $c_3(Y)=-64$\hidewidth&& && &\cr
\tablerule\tablerule
&& && && && && && && &&
&&\omit\hidewidth $|K+H|$ defines a birational\hidewidth
&&\omit\hidewidth unirational\hidewidth&&[Is]&\cr
&&\omit\hidewidth$X_{27}$\hidewidth&&\omit\hidewidth 13\hidewidth
&&\omit\hidewidth 19\hidewidth&&0&&1&&\omit\hidewidth 11\hidewidth
&&\omit\hidewidth $-\infty$\hidewidth
&&\omit\hidewidth$X_{27}\link 4 5 X_{6}$\hidewidth
&&\omit\hidewidth map onto $\openG(1,5)\cap\openP^9$;\hidewidth&&
\omit\hidewidth not rational\hidewidth&&[Ch3]&\cr
&& && && && && && && && &&
\omit\hidewidth birat. to cubic 3-fold in $\openP^4$ [Fa2]\hidewidth&& && &\cr
\tablerule\tablerule
&& && && && && && && && &&
\omit\hidewidth $|K+H|$ defines a birational\hidewidth
&&\omit\hidewidth birational\hidewidth&& &\cr
&&\omit\hidewidth$X_{28}$\hidewidth&&\omit\hidewidth 14\hidewidth
&&\omit\hidewidth 22\hidewidth&&1&&0&&\omit\hidewidth 14\hidewidth&&$0$
&&\omit\hidewidth$X_{28}\link 5 5 X_{17}$\hidewidth
&&\omit\hidewidth map onto $\openG(1,6)\cap\openP^{13}$;\hidewidth
&&\omit\hidewidth Calabi-Yau\hidewidth&&[Ch3]&\cr
&& && && && && && && && && && && &\cr
\tablerule\tablerule
&& && && && && && && && &&
\omit\hidewidth $K=K_1+K_2$, $H^2K_1=10$;\hidewidth&&\omit\hidewidth
birational\hidewidth&& &\cr
&&\omit\hidewidth$X_{29}$\hidewidth&&\omit\hidewidth 17\hidewidth
&&\omit\hidewidth 32\hidewidth&&1&&0&&\omit\hidewidth 24\hidewidth&&$0$&&
\omit\hidewidth $X_{29}\link 5 5 $\hidewidth&&\omit\hidewidth $|H+K+K_1|$
birat. onto a c.i.\hidewidth&&
\omit\hidewidth Calabi-Yau\hidewidth&& &\cr
&& && && && && && &&
&&\omit\hidewidth$\displaystyle{X_1\un\un_{i=1}^5\openP^3}$&&
\omit\hidewidth $(1,1,1,1,1)^2$ in $\displaystyle{\kreuz_{i=1}^5
\openP^1\subset\openP^{31}}$\hidewidth
&&(elliptic)&& &\cr
\tablerule\tablerule
&& && && && && && && &&
&&\omit\hidewidth log-general type\hidewidth&&\omit\hidewidth not\hidewidth&&
&\cr
&&\omit\hidewidth$X_{30}$\hidewidth&&\omit\hidewidth 18\hidewidth
&&\omit\hidewidth 35\hidewidth&&0&&2&&\omit\hidewidth 26\hidewidth
&&\omit\hidewidth $-\infty$\hidewidth&& && &&
\omit\hidewidth rational\hidewidth&& &\cr
&& && && && && && && && && && && &\cr
\tablerule\tablerule
\hfil\cr}}}\endgather
$$
\par\noindent
{\bf Remark 7.4.}\quad (i) The classification of smooth 3-folds in
$\openP^5$ is complete up to degree 11 and  almost complete
in degree 12. \par\noindent
(ii) Some of the information in Table 7.3 is new. It can be obtained along the
lines of section 5 (compare Example 5.8).\par\noindent
(iii) In order to construct $X_{21}$ via a liaison $X_{21}\link 5 5 X_{15}\cup
X_1$,
one has to choose $X_{15}$ and $X_1$ in a special position.\quad\qed
\par\medskip
%
%
%

\enddocument